\documentclass[11pt]{amsart}

\makeatletter
\usepackage{amssymb}
\usepackage{latexsym}
\usepackage{amsbsy}
\usepackage{amsfonts}

\def\marginpar#1{\ignorespaces}

\newtheorem{theorem}[equation]{Theorem}
\newtheorem{proposition}[equation]{Proposition}
\newtheorem{lemma}[equation]{Lemma}
\newtheorem{corollary}[equation]{Corollary}
\newtheorem{definition}[equation]{Definition}

\theoremstyle{definition}
\newtheorem{remark}[equation]{Remark}

\newtheorem{example}[equation]{Example}

\numberwithin{equation}{section}

\def\AArm{\fam0 \rm}%
\newdimen\AAdi%
\newbox\AAbo%
\def\AAk#1#2{\setbox\AAbo=\hbox{#2}\AAdi=\wd\AAbo\kern#1\AAdi{}}%

\newcommand{\BBone}{{\ensuremath{{\AArm 1\AAk{-.8}{I}I}}}}

\def\eqref#1{(\ref{#1})}
\def\eqlabel#1{\def\@currentlabel{#1}}

\def\formula#1{\def\@tempa{#1}\let\@tempb\theequation\def\theequation{%
\hbox{#1}}\def\@currentlabel{(\theequation)}$$}
\def\endformula{\leqno\hbox{(\@tempa)}$$\@ignoretrue\let\theequation\@tempb}

\def\given{\hskip5\p@\relax\vrule\@width.4\p@\hskip5\p@\relax}

\newcommand{\open}[1]{%
\par\normalfont\topsep6\p@\@plus6\p@\trivlist\item[\hskip\labelsep\itshape#1%
\@addpunct{.}]\ignorespaces}

\DeclareRobustCommand{\close}[1]{%
  \ifmmode 
  \else \leavevmode\unskip\penalty9999 \hbox{}\nobreak\hfill
  \fi
  \quad\hbox{$#1$}}

\newlength{\toskip}

\def \R {{\mathbb R}}
\def \Q {{\mathbb Q}}

\def \T {{\mathbb T}}
\def \P {{\mathbb P}}
\def \E {{\mathbb E}}
\def \N {{\mathbb N}}

\def \L {{\mathbb L}}

\def \ep {\varepsilon}
\def \phi {\varphi}
\def \s {\sigma}
\def \a {\alpha}
\def \b {\beta}
\makeatother

\begin{document}
\date{\today}

\title[DEVIATIONS BOUNDS ...]{DEVIATIONS BOUNDS AND CONDITIONAL PRINCIPLES FOR THIN SETS.}

 \author[P. Cattiaux and N. Gozlan]{\quad {Patrick} Cattiaux \, and \, {Nathael} Gozlan }

\address{{\bf {Patrick} CATTIAUX},\\ Ecole Polytechnique, CMAP, F- 91128 Palaiseau cedex,
CNRS 756\\ and Universit\'e Paris X Nanterre, \'equipe MODAL'X, UFR SEGMI\\ 200 avenue de la
R\'epublique, F- 92001 Nanterre, Cedex.} \email{cattiaux@cmapx.polytechnique.fr}
\address{{\bf {Nathael} GOZLAN},\\ Universit\'e Paris X Nanterre, \'equipe MODAL'X, UFR SEGMI\\ 200 avenue de la
R\'epublique, F- 92001 Nanterre, Cedex.}
\email{nathael.gozlan@u-paris10.fr}

\maketitle

 \begin{center}
 \textsc{Ecole Polytechnique \quad and \quad Universit\'e Paris X}
 \end{center}

\begin{abstract}
The aim of this paper is to use non asymptotic bounds for the probability of rare events in the
Sanov theorem, in order to study the asymptotics in conditional limit theorems (Gibbs conditioning
principle for thin sets). Applications to stochastic mechanics or calibration problems for
diffusion processes are discussed.
\end{abstract}
\bigskip

\bigskip

\section{\bf Introduction}\label{I}

Let $X_1,X_2,\ldots$ be i.i.d. random variables taking their values in some metrizable space $(E,d)$.
Set $M_n=\frac 1n \, \sum_{i=1}^n \, X_i$ the empirical mean (assuming here that $E$ is a vector
space) and $L_n=\frac 1n \, \sum_{i=1}^n \, \delta_{X_i}$ the empirical measure. In recent years new
efforts have been made in order to understand the asymptotic behavior of laws conditioned by some
rare or super-rare event.
\smallskip

The celebrated Gibbs conditioning principle is the corresponding meta principle for the empirical
measure, namely $$\lim_{n\rightarrow \, +\infty} \, \P^{\otimes n}((X_1,\ldots,X_k)\in B \, / \, L_n\in
A)= (\mu^*)^{\otimes k}(B) \, ,$$ where $\mu^*$ minimizes the relative entropy $H(\mu^*\mid\mu)$ among
the elements in $A$. When $A$ is thin (i.e. $\P^{\otimes n}(L_n\in A)=0$), such a statement is
meaningless, so one can either try to look at regular desintegration (the so called ``thin shell''
case) or look at some enlargement of $A$. The first idea is also meaningless in general (see however
the work by Diaconis and Freedman \cite{DF88}). Therefore we shall focus on the second one.

An enlargement $A_{\varepsilon}$ is then a non thin set containing $A$, and the previous statement
becomes a double limit one i.e. $$\lim_{\varepsilon \, \rightarrow \,  0} \quad \lim_{n\rightarrow
\, +\infty} \, .$$ Precise hypotheses are known for this meta principle (``thick shell'' case) to
become a rigorous result, and refinements (namely one can choose some increasing $k(n)$) are known
(see e.g. \cite{DZ96} and the references therein). One possible way to prove this result is to
identify relative entropy with the rate function in the Large Deviations Principle for empirical
measures (Sanov's theorem). In this paper we will introduce an intermediate ``approximate thin
shell'' case, i.e. we will look at the case when the enlargement size depends on $n$, i.e.
$\varepsilon_n \, \rightarrow \, 0$. We shall also discuss in details one case of ``super-thin''
set, i.e. when relative entropy is infinite for any element in $A$.

Of course since we are considering conditional probabilities, we are led to get both lower and
upper non asymptotic estimates for the probability of rare events.
\smallskip

The paper is organized as follows.
\smallskip

In Section \ref{Rappels} we shall introduce the notations and recall some results we shall use
repeatedly. Then we give the main general result (Theorem \ref{2.7}).

When $A$ is some closed subspace (i.e. defined thanks to linear constraints), our program can be
carried out by directly using well known inequalities for the sum of independent variables. This
will be explained in Section \ref{Moment}.

The more general case of a general convex constraint is studied in Section \ref{Convex}. In the
compact case upper estimates are well known and lower estimates will be derived thanks to a result
by Deuschel and Stroock. In both cases on has to compute the metric entropy (i.e. the number of
small balls needed to cover $A$) for some metric compatible with the convergence of measures. The
extension to non compact convex constraint is done by choosing an adequate rich enough compact
subset.

Section \ref{Exemples} is devoted to some examples, first in a finite dimensional space. We next
show that the Schr\"{o}dinger bridges and the Nelson processes studied in Stochastic Mechanics, are
natural ``limiting processes'' for constraints of marginal type.

Section \ref{Calib} is devoted to the study of a super-thin example corresponding to the well known
problem of volatility calibration in Mathematical Finance. Our aim is to give a rigorous status to
the ``Relative Entropy Minimization method'' introduced in \cite{AFHS}. The problem here is to
choose the diffusion coefficient (volatility) of a diffusion process with a given drift (risk
neutral drift), knowing some final moments of the diffusion process. Of course all the possible
choices are mutually singular so that the constraint set $A$ does not contain any measure with
finite relative entropy, i.e. is super-thin. We shall show that under some conditions, the method
by Avellaneda et altri \cite{AFHS} enters our framework, hence furnishes the natural candidate from
a statistical point of view (we shall not discuss any kind of financial related aspects).
\medskip

Another famous example of super-thin set is furnished by Statistical Mechanics, namely: are the
Gibbs measures associated to some Hamiltonian the limiting measures of some conditional law of large
numbers ? The positive answer gives an interpretation of the famous Equivalence of Ensembles
principle (see \cite{SZ91,DSZ}). It should be interesting to relate the Gibbs variational principle
 as in \cite{DSZ} to the above Gibbs conditioning principle. This is not done here.
\medskip

{\bf Acknowledgements.} We want to warmly acknowledge Christian L\'eonard for so many animated
conversations on Large Deviations Problems, and for indicating to us various references on the
topic.
\bigskip

\section{\bf Notation and first basic results.}\label{Rappels}

Throughout the paper $(E,d)$ will be a Polish space. $M_1(E)$ (resp. $M(E)$) will denote the set of
Probability measures (resp. bounded signed measures) on $E$ equipped with its Borel $\sigma$-field.
$M_1(E)$ is equipped with the narrow topology (convergence in law) and its natural Borel
$\sigma$-field.

In the sequel, we will consider a sequence $X_1,X_2,\ldots$ of i.i.d. $E$ valued random variables. The common law of the $X_i$'s will be denoted by $\alpha$ and their empirical measure by $L_n=\frac 1n \, \sum_{i=1}^{n} \, \delta_{X_i}$. 

\medskip

Our aim is to study the asymptotic behavior of the conditional law
\begin{equation}\label{2.1}
\alpha_{A,k}^n(B)=\P_{\alpha}^{\otimes n}\Big((X_1,\ldots,X_k)\in B \, / \, L_n\in A_n\Big)
\end{equation}
for some $A_n$ going to some thin set $A$ when $n$ goes to $\infty$.
\smallskip

The first tool we need is relative entropy. Recall that for $\b$ and $\gamma$ in $M_1(E)$, the
relative entropy $H(\b\mid\gamma)$ is defined by the two equivalent formulas
\begin{enumerate}\stepcounter{equation}\eqlabel{\theequation}\label{2.2}
\item[(\theequation.1)] \quad $H(\b\mid\gamma)=\int \, \log  \big(\frac{d\b}{d\gamma}\big) \,
d\b$ , if this quantity is well defined and finite, $+\infty$ otherwise,
\item[(\theequation.2)] \quad $H(\b\mid\gamma)=\sup \, \{\int \phi d\b-\log \int e^\phi \,
d\gamma \, , \, \phi \in C_b(E)\}$ .
\end{enumerate}
If $B$ is a measurable set of $M_1(E)$ we will write
\begin{equation}\label{2.3}
H(B\mid\gamma)=\inf \, \{H(\b\mid\gamma) \, , \, \b \in B\} \, .
\end{equation}
\smallskip

The celebrated Sanov's theorem tells that for any measurable set $B$ $$- \, H(\overset{\circ}
B\mid\alpha) \, \leq \, \liminf_{n\to \infty} \, \frac 1n \, \log \, \mathbb P(L_n \in B) \, \leq \,
\limsup_{n\to \infty} \, \frac 1n \, \log \, \mathbb P(L_n \in B) \, \leq \, - \, H(\overline B\mid\alpha)
\, ,$$ where the interior $\overset{\circ} B$ and the closure $\overline B$ of $B$ are for the narrow
topology.

Recall that one can reinforce the previous topology by considering the $G$-topology induced by some
subset $G$ of measurable functions containing all the bounded measurable functions. In particular if
$\alpha$ satisfies the strong Cram\'er assumption i.e $\forall g \in G \, , \, \forall t
> 0 \, $,
\begin{equation}\label{eq:Cramer}
\int \, e^{tg} \, d\alpha \, < \, +\infty \, ,
\end{equation}
the previous result is still true for the $G$-topology (see \cite{EiS} thm. 1.7). When $G$ is
exactly the set of measurable and bounded functions, the $G$-topology is usually called
$\tau$-topology.
\smallskip

It is thus particularly interesting to have some information on the possible Arginf in \eqref{2.3}.
The result below is collecting some known facts:

\begin{theorem}\label{2.5}
Let $C$ be a measurable convex subset of $M_1(E)$ such that $H(C\mid\alpha) < +\infty$. There exists an
unique probability measure $\alpha^*$ such that any sequence $\nu_n$ of $C$ such that $\lim_{n \to
+\infty} \, H(\nu_n\mid\alpha) = H(C\mid\alpha)$, converges in total variation distance to $\alpha^*$.

This probability measure (we shall call the generalized $I$- projection of $\alpha$ on $C$) is
characterized by the following Pythagoras inequality $H(\nu\mid\alpha) \, \leq \, H(\alpha^*\mid\alpha) +
H(\nu\mid\alpha^*)$ for all $\nu \in C$.

If $\alpha^*$ belongs to $C$ we shall call it the $I$- projection (non generalized). In particular
the $I$- projection on a total variation closed convex subset such that $H(C\mid\alpha) < +\infty$
always exists.

Finally if $\alpha$ satisfies the strong Cramer assumption \eqref{eq:Cramer} one can replace total
variation closed by $G$-closed in the previous statement.
\end{theorem}

All these results can be found in \cite{Csi75,EiS} (see \cite{G}, chap. II for more details).
\smallskip
\pagebreak

Before to state our first results on thin constraints, we recall the known results on thick ones.

\begin{theorem}\eqlabel{\theequation}\label{2.6}
\begin{enumerate}
\item[] \item[(\theequation.1)] (see \cite{Csi84}). If $C$ is convex, closed for the $\tau$-topology
and satisfies $H(C\mid\alpha)=H(\overset{\circ} C\mid\alpha)<+\infty$ then $\alpha_{C,k}^n$ defined in
\eqref{2.1} is well defined for $n$ large enough, and converges (when $n$ goes to $\infty$) in
relative entropy to $\alpha^{* \, \otimes k}$, where $\alpha^*$ is the $I$- projection of $\alpha$
on $C$. \item[(\theequation.2)] (see \cite{SZ91}). If $A$ is a measurable subset such that $H(\bar
A\mid\alpha)=H(\overset{\circ} A\mid\alpha)<+\infty$, and if there exists an unique $\alpha^* \in \bar A$
such that $H(\bar A\mid\alpha)=H(\alpha^*\mid\alpha)$, then $\alpha_{A,k}^n$ again converges to $\alpha^{*
\, \otimes k}$ but for the narrow convergence.
\end{enumerate}
\end{theorem}

When $H(\overset{\circ} A\mid\alpha)= +\infty$ (in particular if $\overset{\circ} A$ is empty) but
$H(A\mid\alpha) < +\infty$ (thin constraints) we have to face some new problems. The strategy is then
to enlarge $A$, considering some nice $A_\varepsilon$, and to consider limits first in $n$, next in
$\varepsilon$. Here we shall consider enlargements depending on $n$. Here is a general result in
this direction.

\begin{theorem}\eqlabel{\thetheorem}\label{2.7}
\quad Let $C_n$ be a non increasing sequence of convex subsets, closed for the $G$-topology. Denote
by $C \, = \, \bigcap_{n=1}^{\infty} \, C_n$. Assume that
\begin{enumerate}
\item[(\thetheorem.1)] $H(C\mid\alpha) \, < \, +\infty$ , \item[(\thetheorem.2)] $\alpha$ has an $I$-
projection $\alpha^*$ on $C$, \item[(\theequation.3)] $\lim_{n \to \infty} \, H(C_n\mid\alpha) \, = \,
H(C\mid\alpha)$ , \item[(\thetheorem.4)] $\liminf_{n \to \infty} \, \frac 1n \, \log \alpha^{\otimes
n}(L_n \in C_n) \, \geq \, - \, H(C\mid\alpha)$ .
\end{enumerate}
Then, for all $k \in \mathbb N^*$, $\alpha_{C_n,k}^n$ converges in total variation distance to
$\alpha^{* \, \otimes k}$.
\end{theorem}

\begin{remark}\label{2.8}
Define $$\mathbb L^a_{\tau}(\alpha) \, = \, \left\{ \, g \textrm{ measurable } : \forall s \in
\mathbb R \, , \, \int \, e^{s | g|} \, d\alpha \, < \, +\infty \, \right\} \, .$$ If $G\subseteq
\mathbb L^a_{\tau}(\alpha)$, we already know (see Theorem \ref{2.5}) that $\alpha^*$ exists as soon
as $H(C\mid \alpha)$ is finite. In addition, since the relative entropy is a good rate function
(according to \cite{EiS} its level sets are compact) (\ref{2.7}.3) is also satisfied. Hence, in
this case, assuming $H(C\mid \alpha)<+\infty$, the only remaining condition to check is
\begin{equation}\label{2.9}
\liminf_{n \to \infty} \, \frac 1n \, \log \alpha^{\otimes n}(L_n \in C_n) \, \geq \, - \,
H(C\mid \alpha) \, .
\end{equation}
\end{remark}
\smallskip

\begin{proof}{\textit {of Theorem \ref{2.7}.}}

Let $\alpha_n^*$ the generalized $I$- projection of $\alpha$ on $C_n$. then
\begin{eqnarray}\label{2.10}
\parallel \alpha_{C_n,k}^n \, - \, \alpha^{* \, \otimes k}\parallel_{TV} & \leq &
\parallel \alpha_{C_n,k}^n \, - \, \alpha_n^{* \, \otimes k}\parallel_{TV} \, + \,
\parallel \alpha_n^{* \, \otimes k} \, - \, \alpha^{* \, \otimes k}\parallel_{TV} \\ & \leq &
\sqrt{2 H\left(\alpha_{C_n,k}^n\mid \alpha_n^{* \, \otimes k}\right)} \, + \, \sqrt{2 H\left(\alpha^{*
\, \otimes k}\mid \alpha_n^{* \, \otimes k}\right)} \nonumber \\ & \leq & \sqrt{2
H\left(\alpha_{C_n,k}^n\mid \alpha_n^{* \, \otimes k}\right)} \, + \, \sqrt{2k H\left(\alpha^*
\mid \alpha_n^* \right)} \nonumber
\end{eqnarray}
where we have used successively the triangle inequality, Pinsker inequality and the additivity of
relative entropy. Since $\alpha^*$ is the $I$- projection of $\alpha$ on $C$, $\alpha^*$ belongs to
$C$ and all $C_n$, so that using Theorem \ref{2.5}, $$H(C\mid \alpha) \, = \, H(\alpha^*\mid \alpha) \, \geq
\, H(\alpha^*\mid \alpha^*_n) \, + \, H(C_n\mid \alpha) \, .$$ Thanks to (\ref{2.7}.3) we thus have $\lim_{n
\to \infty} \, H(\alpha^*\mid \alpha^*_n) \, = \, 0$ .

To finish the proof (according to \eqref{2.10}) it thus remains to show that $\lim_{n \to \infty}
\, H\left(\alpha_{C_n,k}^n\mid \alpha_n^{* \, \otimes k}\right) \, = \, 0$ . But thanks to
(\ref{2.7}.4), for $n$ large enough, $\alpha^{\otimes n}(L_n \in C_n) > 0$, so that we may apply
Lemma \ref{2.11} below with $A=C_n$. It yields
\begin{eqnarray*}
H\left(\alpha_{C_n,k}^n\mid \alpha_n^{* \, \otimes k}\right) & \leq & - \, \frac kn \, \log
\left(\alpha^{\otimes n}(L_n \in C_n) \, e^{n H(C_n\mid \alpha)} \right) \\ & \leq & - \, \frac kn \,
\log \left(\alpha^{\otimes n}(L_n \in C_n) \, e^{n H(C\mid \alpha)} \right) \, + \, k \,
\left(H(C\mid \alpha)-H(C_n\mid \alpha)\right) \, .
\end{eqnarray*}

According to (\ref{2.7}.4) the first term in the right hand side sum has a non positive $\limsup$,
while the second term goes to $0$ thanks to (\ref{2.7}.3). Since the left hand side is nonnegative the
result follows.
\end{proof}

We now recall the key Lemma due to Csiszar (\cite{Csi84}) we have just used :

\begin{lemma}\label{2.11}
Let $A$ be a convex $G$-closed subset, such that $H(A\mid \alpha)<+\infty$. Denote by $\alpha^*$ the
generalized $I$- projection of $\alpha$ on $A$. Then if $\alpha^{\otimes n}(L_n \in A) > 0$, for all
$k \in \mathbb N^*$ , $$H\left(\alpha_{A,k}^n\mid \alpha^{* \, \otimes k}\right) \, \leq \, - \,
\frac{1}{[n/k]} \, \log \left(\alpha^{\otimes n}(L_n \in A) \, e^{n H(A\mid \alpha)} \right) \, .$$
\end{lemma}
\medskip

Under some additional assumption one can improve the convergence in Theorem \ref{2.7}. Introduce
the usual Orlicz space $$\mathbb L_{\tau}(\alpha) \, = \, \left\{ \, g \textrm{ measurable } :
\exists s \in \mathbb R \, , \, \int \, e^{s |g|} \, d\alpha \, < \, +\infty \, \right\} \, .$$
Note the difference with $\L^a_{\tau}$ (for which $\exists$ is replaced by $\forall$). We equip
$\L_{\tau}$ with the Luxemburg norm $$\parallel g\parallel_{\tau} \, = \, \inf \, \{s>0 \, , \,
\int \tau(g/s) \, d\alpha \, \leq \, 1\} \textrm{ where } \tau(u)=e^{|u|} - |u| - 1 \, .$$ It is
well known that the dual space of $\L_{\tau}(\alpha)$ contains the set of probability
measures $\nu$ such that $H(\nu\mid \alpha)<+\infty$. We equip this dual space with the dual norm
$\parallel \, \parallel^*_{\tau}$.

\begin{proposition}\label{2.12}
In addition to all the assumptions in Theorem \ref{2.7}, assume the following: the densities $h_n
\, = \, \frac{d\alpha_n^*}{d\alpha}$ ($\alpha_n^*$ being the generalized $I$- projection of
$\alpha$ on $C_n$) define a bounded sequence in $\L^p(\alpha)$ for some $p>1$. Then
$$\lim_{n \to +\infty} \parallel \alpha_{C_n}^n \, - \, \alpha^*\parallel^*_{\tau} \, = \, 0 \, .$$
\end{proposition}

\begin{proof}
The proof is exactly the same as the one of Theorem \ref{2.7} with $k=1$, just replacing $\parallel
\,
\parallel_{TV}$ by $\parallel \, \parallel^*_{\tau}$ in the first line of \eqref{2.10}, and then
replacing Pinsker inequality by the following one, available for $\nu_i$'s such that
$H(\nu_i\mid \alpha)<+\infty$,
\begin{equation}\label{2.13}
\parallel \nu_1 \, - \, \nu_2 \parallel^*_{\tau} \, \leq \, q C \, (1+\log(4^{1/q} \parallel
\frac{d\nu_2}{d\alpha}\parallel_p)) \, \left(H(\nu_1\mid \nu_2)+\sqrt{H(\nu_1\mid \nu_2)}\right) \, ,
\end{equation}
where $q=p/(p-1)$, $\nu_2$ being $\alpha_n^*$ and $\nu_1$ being either $\alpha_{C_n}^n$ or
$\alpha^*$.

In order to prove \eqref{2.13} we first recall the weighted Pinsker inequality recently shown by
Bolley and Villani \cite{BV05} (also see \cite{G} for another approach) : there exists some $C$
such that for all nonnegative $f$ and all $\delta > 0$ , $$\parallel f\nu_1 \, - \, f
\nu_2\parallel_{TV} \, \leq \, (C/\delta) \, \left(1 + \log \, \int e^{\delta f} d\nu_2\right) \,
\left(H(\nu_1\mid \nu_2)+\sqrt{H(\nu_1\mid \nu_2)}\right) \, .$$ For a $f$ such that $\parallel
f\parallel_{\tau} \leq 1$ and $\delta=1/q$ it thus holds, first $\int e^{|f|} d\alpha \leq 4$, then
thanks to H\"{o}lder's inequality $\int e^{\delta |f|} d\nu_2 \, \leq \,4^{\delta} \parallel
\frac{d\nu_2}{d\alpha}\parallel_p$. \eqref{2.13} immediately follows.
\end{proof}

\section{\bf $F$ moment constraints.}\label{Moment}

In this section $G=\mathbb L_{\tau}(\alpha)$ and we consider constraints $C$ in the form $$C \, =
\, \left\{\nu\in M_1(E) \, , \, \int \, F \, d\nu \, \in \, K\right\}$$ where $F$ is a measurable $B$ valued
map ($(B,\parallel\,.\, \parallel)$ being a separable Banach space equipped with its cylindrical
$\sigma$-field) where $\int F d\nu$ denotes the Bochner integral and $K$ is a closed convex set of
$B$. We denote by
$$\forall \lambda \in B',\quad Z_F(\lambda)=\int_E \, \exp(\langle\lambda, F(x)\rangle) \, \alpha(dx) \, ,\qquad \Lambda_F(\lambda)=\log Z_F(\lambda) $$ the
Laplace transform and moment generating function of $F$.

We always assume that
\begin{itemize}
\item  $\parallel F
\parallel \, \in \, \mathbb L_{\tau}(\alpha)$, \item $dom\, \Lambda_F \, = \, \{\lambda \in B' \, ,
\, \Lambda_F(\lambda)<+\infty\}$ ($B'$ being the dual space of $B$) is a non empty open set of
$B'$.
\end{itemize}

The enlargement $C_n$ is defined similarly $$C_n \, = \, \left\{\nu\in M_1(E) \, , \, \int \, F \, d\nu
\, \in \, K^{\varepsilon_n} \right\}$$ for $K^{\varepsilon_n} \, = \, \{x\in B \, , \, d(x,C_n)\leq
\varepsilon_n\}$.
\medskip

What we have to do is to check all the assumptions of Theorem \ref{2.7}. But the situation here is
particular since the condition $L_n \in A$ reduces to $\sum_{i=1}^n \, F(X_i) \, \in \, A$. Thanks
to the next Lemma \ref{3.1} assumption (\ref{2.7}.4) reduces to well known estimates:

\begin{lemma}\label{3.1}
Assume that the $I$- projection $\alpha^*$ of $\alpha$ on $C$ exists and can be written $\alpha^*
\, = \, \frac{e^{\langle\lambda^*,F\rangle}}{Z_F(\lambda^*)} \, \, \alpha$ for some $\lambda^* \in B'$ . Then
for all $\varepsilon >0$, $$\frac 1n \, \log \left(\alpha^{\otimes n}(L_n \in C_\varepsilon) \,
e^{n H(\alpha^*\mid \alpha)}\right) \, \geq \, \frac 1n \, \log \mathbb P \left(\left\| \frac 1n \,
\sum_{i=1}^n \, F(Y_i) \, - \, \int  F \, d\alpha^*\right\| \leq \varepsilon\right) \, - \,
\parallel \lambda^*\parallel \varepsilon \, ,$$ where the $Y_i$'s are i.i.d. random variables with
common law $\alpha^*$.
\end{lemma}

\begin{proof}
The proof uses the standard centering method in large deviations theory. Denote by $L_n^x=\frac 1n
\sum_{i=1}^n \delta_{x_i}$ the empirical measure of $x=(x_1,\ldots,x_n)$. Then
\begin{eqnarray*}
\alpha^{\otimes n}(L_n \in C_\varepsilon) & = & \int \BBone_{C_\varepsilon}(L_n^x) \, \prod_{i=1}^n
\frac{d\alpha}{d\alpha^*}(x_i) \, d\alpha^{* \otimes n}(x) \\ & = & \int
\BBone_{C_\varepsilon}(L_n^x) \, \exp\left(-n\left\langle L_n^x,\log \frac{d\alpha}{d\alpha^*}\right\rangle\right) \,
d\alpha^{* \otimes n}(x) \\ & = & e^{-n H(\alpha^*\mid \alpha)} \, \int \BBone_{C_\varepsilon}(L_n^x)
\, \exp\left(-n\left\langle L_n^x-\alpha^*,\log \frac{d\alpha}{d\alpha^*}\right\rangle \right) \, d\alpha^{* \otimes n}(x)
\\ & = & e^{-n H(\alpha^*\mid \alpha)} \, \int \BBone_{C_\varepsilon}(L_n^x)
\, \exp\left(-n\left\langle \lambda^*,\frac 1n \, \sum_{i=1}^n F(x_i)  -  \int  F d\alpha^*\right\rangle\right) \,
d\alpha^{* \otimes n}(x)
\end{eqnarray*}
Now we may replace $C_\varepsilon$ by its subset $$\widetilde C_\varepsilon = \left\{\nu \in M_1(E) \, , \,
\int \parallel F\parallel d\nu <+\infty \textrm{ and } \left\| \int F d\nu - \int F
d\alpha^*\right\| \, \leq \, \varepsilon\right\}$$ and obtain
\begin{eqnarray*}
\alpha^{\otimes n}(L_n \in C_\varepsilon)e^{n H(\alpha^*\mid \alpha)}&\geq&\int \BBone_{\widetilde
C_\varepsilon}(L_n^x) e^{-n\left\langle\lambda^*,\frac 1n \, \sum_{i=1}^n F(x_i)  -  \int  F d\alpha^*\right\rangle}
d\alpha^{* \otimes n}(x)\\ &\geq& e^{-n\parallel \lambda^*\parallel \varepsilon} \, \int
\BBone_{\widetilde C_\varepsilon}(L_n^x) \, d\alpha^{* \otimes n}(x)
\end{eqnarray*}
that completes the proof.
\end{proof}

The next Lemma \ref{3.2} is well known in convex analysis. For a complete proof the reader is
referred to \cite{G} Lemma II.39,

\begin{lemma}\label{3.2}
Under our hypotheses on $F$ and $dom\,\Lambda_F$, if in addition the function
$$H(\lambda)=\Lambda_F(\lambda)-\inf_{y\in K}\langle\lambda,y\rangle$$ achieves its minimum at (at least one)
$\lambda^*$, then $H(C\mid \alpha) \, = \, \sup_{\lambda \in B'} \, \{\inf_{y \in K} \langle\lambda,y\rangle -
\Lambda_F(\lambda)\}$ and the $I$- projection $\alpha^*$ of $\alpha$ on $C$ exists and can be
written $\alpha^* \, = \, \frac{e^{\langle\lambda^*,F\rangle}}{Z_F(\lambda^*)} \, \, \alpha$.

In the sequel we shall denote (H-K) the additional assumption on $H$.  In particular if $K=\{x_0\}$
with $x_0=\nabla \Lambda_F(\lambda_0)$ (H-K) is satisfied.
\end{lemma}

Before to state our first general result let us recall some definition

\begin{definition}\label{3.3}
$B$ is of type 2 if there exists some $a>0$ such that for all sequence $Z_i$ of $\mathbb L^2$
i.i.d. random variables with zero mean and variance equal to 1, the following holds $$\mathbb E
\left(\left\| \sum_{i=1}^n Z_i\right\|^2\right) \leq a \, \sum_{i=1}^n \, \mathbb E(\parallel
Z_i \parallel^2) \, .$$ In particular an Hilbert space is of type 2.
\end{definition}

We arrive at

\begin{theorem}\label{3.4}
In addition to our hypotheses on $F$ and $dom \Lambda_F$, assume that $B$ is of type 2 and that
(H-K) is satisfied. If $\varepsilon_n \, > \, \frac{c}{\sqrt n}$ with $c=\sqrt{a \,
Var_{\alpha^*}(F)}$, then $\alpha_{C_n,k}^n$ converges to $\alpha^{* \otimes k}$ in total variation
distance when $n \to \infty$.
\end{theorem}

\begin{proof}
(\ref{2.7}.1) and (\ref{2.7}.2) are satisfied with our hypotheses, according to Lemma \ref{3.2}.

In order to prove (\ref{2.7}.3) introduce the function $H_n$ similar to $H$ in Lemma \ref{3.2}
replacing $K$ by $K^{\varepsilon_n}$. Of course $$\inf H \, \leq \, \inf H_n \, \leq \,
H_n(\lambda^*) \, \xrightarrow[n\rightarrow+\infty]{} \, \inf H \, = \, H(\lambda^*)$$ since $H_n$ converges to $H$
pointwise on the domain of $H$. We already know that $\inf H = - H(C\mid \alpha)$. It is thus enough to
prove that $\inf H_n = - H(C_n\mid \alpha)$. But this is a consequence of Csiszar results (\cite{Csi75}
thm 3.3 and \cite{Csi84} thm 2 and 3, also see \cite{G} thm II.41 for another proof) since the
intersection of the interior of  $K^{\varepsilon_n}$ and the convex hull of the support of the
image measure $F^{-1} \alpha$ is non empty.

Finally in order to prove (\ref{2.7}.4), according to Lemma \ref{3.1} it is enough to check that
$$\lim_{n \to \infty} \, \frac 1n \, \log \mathbb P \left(\left\| \frac 1n \, \sum_{i=1}^n \, F(Y_i) \, - \, \int  F
\, d\alpha^*\right\| \leq \varepsilon_n\right) \, = \, 0 \, .$$ To this end recall the following
theorem of Yurinskii

\begin{theorem}\label{3.5} \textrm{(Yurinskii, \cite{Yur76} theorem 2.1)}. If $Z_i$ is a $B$ valued
sequence of centered independent variables such that there exist $b$ and $M$ both positive, with
$$\forall i \in \mathbb N^* \, , \, \forall k \geq 2, \quad \mathbb E(\parallel Z_i\parallel^k) \,
\leq \, \frac{k!}{2} \, b^2 M^{k-2} \, ,$$ then denoting $S_n = \sum_{i=1}^n \, Z_i$ it holds
$$\forall t>0, \, \mathbb P(\parallel S_n \parallel \geq \mathbb E(\parallel S_n\parallel) + nt) \, \leq \,
\exp \, \left( - \, \frac 18 \, \frac{nt^2}{b^2+tM}\right) \, .$$
\end{theorem}

We may apply Theorem \ref{3.5} with $Z_i=F(Y_i)-\int F d\alpha^*$, $M=\parallel F - \int F
d\alpha^*\parallel_{\L_\tau(\alpha^*)}$ and $b=\sqrt{ 2} M$ as soon as $F \in \L_\tau(\alpha^*)$.
Indeed since $B$ is of type 2, $\mathbb E(\parallel S_n\parallel)\leq \sqrt{\mathbb E(\parallel
S_n\parallel^2)} \leq \sqrt{a n} \sigma$ with $\sigma=\sqrt{\mathbb E(\parallel Z_1\parallel^2)}$.
It follows $$\mathbb P \left(\left\| \frac 1n \, \sum_{i=1}^n \, F(Y_i) \, - \, \int  F \,
d\alpha^*\right\| \leq \frac{\sigma \sqrt a}{\sqrt n}(1+ t)\right) \, \geq \, 1 \, - \, \exp \,
\left(- \, \frac 18 \, \frac{a \sigma^2 t^2}{2M^2+tM}\right) \, ,$$ and the result provided
$\varepsilon_n \sqrt n > \sigma \sqrt a$.

It remains to prove that $F \in \L_\tau(\alpha^*)$. But thanks to the representation of $\alpha^*$
obtained in Lemma \ref{3.2}
\begin{eqnarray*}
\int e^{t \parallel F\parallel} d\alpha^* &=& \frac{1}{Z_F(\lambda^*)} \, \int e^{t \parallel
F\parallel} \, e^{\langle\lambda^*,F\rangle} \, d\alpha \\ &\leq & \frac{1}{Z_F(\lambda^*)} \, \left(\int e^{tq
\parallel F\parallel} d\alpha\right)^{\frac 1q} \, \left(\int e^{\langle p\lambda^*,F\rangle}
d\alpha\right)^{\frac 1p} \, .
\end{eqnarray*}
Since $dom\,\Lambda_F$ is a non empty open set containing $\lambda^*$, there exists some $p>1$ such
that $p \lambda^* \in dom\,\Lambda_F$, and the result follows for $t$ small enough since $F \in
\L_\tau(\alpha)$.
\end{proof}

\begin{remark}\label{3.6}
Note that if for instance $F$ is bounded everything in Theorem \ref{3.4} can be explicitly
described with the only parameter $n$. However (unfortunately) we do not know any explicit bound
for the speed of convergence of $\alpha_{C_n,k}^n$, because we do not know in general how to
evaluate $H(C_n\mid \alpha)-H(C\mid \alpha)$. Hence from a practical point of view, if we know how to
enlarge $C$, we do not know when a possible algorithm has to be stopped.
\end{remark}

It is natural to ask whether $\varepsilon_n \approx 1/\sqrt n$ is the optimal order for the
enlargement or not. In one dimension the answer is negative as we shall see below

\begin{theorem}\label{3.7}
If $B=\R$ the conclusion of Theorem \ref{3.4} remains true for $\varepsilon_n > c/n$ for some $c$
large enough.
\end{theorem}
\begin{proof}
We shall just replace Yurinskii's estimate by Berry-Eessen bound. Indeed Berry-Eessen theorem tells
us that $$\mathbb P \left(\left| \frac 1n \, \sum_{i=1}^n \, F(Y_i) \, - \, \int  F \,
d\alpha^*\right| \leq \varepsilon_n\right) \, \geq \, \Phi\left(\frac{\varepsilon_n \sqrt
n}{\sigma}\right) \, - \, \Phi\left(- \, \frac{\varepsilon_n \sqrt n}{\sigma}\right) \, - \, 20 \,
\frac{\kappa}{\sigma^3 \sqrt n}$$ where $\Phi(u)=\int_{-\infty}^u \, e^{-s^2/2} ds/\sqrt{2\pi}$,
$\sigma^2=Var_{\alpha^*} F$ and $\kappa$ is the $\alpha^*$'s moment of order 3 of $F-\int F
d\alpha^*$ . It easily follows that $$\mathbb P \left(\left| \frac 1n \, \sum_{i=1}^n \, F(Y_i)
\, - \, \int  F \, d\alpha^*\right| \leq \varepsilon_n\right) \, \geq \, \frac{2}{\sqrt n} \,
\left(\frac{n \varepsilon_n}{\sqrt{2\pi}} \, e^{-n\varepsilon_n^2/2\sigma^2} \, - \,
10(\kappa/\sigma^3)\right) \, = \, \theta_n \, .$$ The requested $1/n \, \log \theta_n \, \to \, 0$
follows with $\varepsilon_n=c/n$ provided $c>10 \sqrt{2\pi} (\kappa/\sigma^3)$.
\end{proof}

Again one may ask about optimality. Actually it is not difficult to build examples with
$\varepsilon_n=c'/n$ for some small $c'$ such that $\mathbb P (L_n\in C_n)=0$ for all $n$. In a
sense this is some proof of optimality. But we do not know how to build examples such that the
previous probability is not zero.
\smallskip

Finally we may improve the convergence, still in the finite dimensional case under slightly more
restrictive assumption.

\begin{theorem}\label{3.8}
In Theorem \ref{3.4} assume that $B=\R^d$ and replace the hypothesis (H-K) by the following : $K \,
\cap \, \overset{\circ} S \, \neq \, \emptyset$ where $S$ is the convex hull of the support of the
image measure $F^{-1} \alpha$. Then $\alpha_{C_n,k}^n$ converges to $\alpha^{* \otimes k}$ both for
the dual norm $\parallel \, \parallel^*_{\tau}$ and in relative entropy.
\end{theorem}

\begin{proof}
The first point is that the new hypothesis is stronger than (H-K). Indeed it is known (see e.g.
\cite{DZ} or  \cite{G} Lemma III.65 for complete proofs) that not only (H-K) holds (as well as
(H-$K^{\varepsilon_n}$) of course), but the minimizers $\lambda^*$ and $\lambda^*_n$ are unique and
$\lambda^*_n \, \to \, \lambda^*$ as $n \to \infty$. Hence $H(C_n\mid \alpha) \, \to \, H(C\mid \alpha)$
too.

Next $\int \left(\frac{d\alpha^*_n}{d\alpha}\right)^p \, d\alpha =
\frac{Z_F(p\lambda^*_n)}{Z_F^p(\lambda^*_n)}$ . Since $\lambda^*_n$ is a bounded (convergent)
sequence, the above quantity can be easily bounded for some $p>1$ (using again the fact that $dom
\Lambda_F$ is an open set). Convergence for the dual norm $\parallel \, \parallel^*_{\tau}$ follows
from Proposition \ref{2.12}.

Finally using exchangeability we have
\begin{eqnarray*}
H(\alpha_{C_n,k}^n\mid \alpha_n^{* \otimes k}) &=& H(\alpha_{C_n,k}^n\mid \alpha^{* \otimes k}) \, + \,
\int \log \frac{d\alpha^{* \otimes k}}{d\alpha_n^{* \otimes k}} \, d\alpha_{C_n}^n \\ &=&
H(\alpha_{C_n,k}^n\mid \alpha^{* \otimes k}) + k H(\alpha^*\mid \alpha^*_n) + k \int \log \frac{d\alpha^*
}{d\alpha_n^*} \, (d\alpha_{C_n}^n - d\alpha^*) \, .
\end{eqnarray*}
We already saw in the proof of Theorem \ref{2.7} that $H(\alpha^*\mid \alpha^*_n)$ and
$H(\alpha_{C_n,k}^n\mid \alpha_n^{* \otimes k})$ go to 0. It remains to prove that $\int \log
\frac{d\alpha^* }{d\alpha_n^*} \, (d\alpha_{C_n}^n - d\alpha^*)$ goes to 0. But $\log
\frac{d\alpha^* }{d\alpha_n^*}=\left\langle\lambda^*_n - \lambda^*,F\right\rangle$ is bounded in $\L_\tau(\alpha)$ for $n$
large enough since $\lambda_n^*$ goes to $\lambda^*$. Hence convergence to 0 of this last term
follows from the convergence for the dual norm $\parallel \, \parallel^*_{\tau}$ we have just
shown.
\end{proof}

\begin{remark}\label{3.9}
In Theorem \ref{3.8} one can also replace Yurinskii's bound by the classical Bernstein inequality
(see e.g. \cite{VW95}). This only improves the constants (see \cite{G} for the details).
\end{remark}

 The results of this Section are satisfactory mainly thanks to Lemma \ref{3.1} and the
very complete literature on sums of independent variables. The situation is of course more
intricate in more delicate situations. We shall study some of them in the next sections.

\section{\bf General convex constraints.}\label{Convex}

We start with the key minimization bound we shall use. The following result is stated in \cite{DS}
Exercise 3.3.23 p76. A complete proof is contained in \cite{Go05} (also see \cite{G}).

\begin{proposition}\label{4.1}
Let $A\subseteq M_1(E)$ be such that $\{x \, , \, L_n^x\in A\}$ is measurable. If $\nu$ is such
that $\nu \ll \alpha$ and $\nu^{\otimes n}(L_n\in A)>0$, then
\begin{eqnarray*}
\frac 1n \, \log \left(\alpha^{\otimes n}(L_n \in A) \, e^{n H(\nu\mid \alpha)}\right) & \geq & - \,
H(\nu\mid \alpha) \, \frac{\nu^{\otimes n}(L_n \in A^c)}{\nu^{\otimes n}(L_n \in A)} \, + \, \frac 1n
\log {\nu^{\otimes n}(L_n \in A)}\\ & & \quad - \, \frac{1}{n e \nu^{\otimes n}(L_n \in A^c)} \, .
\end{eqnarray*}
\end{proposition}

\begin{corollary}\label{4.2}
If (\ref{2.7}.1,2,3) are all satisfied, (\ref{2.7}.4) holds as soon as $$\lim_{n \to +\infty} \,
\alpha^{* \otimes n}(L_n \in C_n) \, = \, 1 \, .$$
\end{corollary}

The proof is an immediate application of Proposition \ref{4.1} with $A=C_n$ and $\nu=\alpha^*$
since $H(C\mid \alpha)=H(\alpha^*\mid \alpha)$.
\smallskip

In the remainder of the section we shall assume that $G=C_b(E)$. According to Remark \ref{2.8}, it
is thus enough to check (\ref{2.7}.1 and 4) in order to apply Theorem \ref{2.7}. In particular if
$H(C\mid \alpha)$ is finite, it just remains to check the condition stated in Corollary \ref{4.2}, by
choosing appropriate enlargements $C_n$. To this end we first recall basic facts on metrics on
probability measures.
\medskip

Recall that the narrow topology on $M_1(E)$ is metrizable. Among admissible metrics we shall
consider two, namely the Prohorov metric $d_P$ and the Fortet-Mourier metric $d_{FM}$.
\begin{proposition}\label{4.3}
For two probability measures $\nu_1$ and $\nu_2$ on $E$ the previous metrics are defined as follows
$$d_P(\nu_1,\nu_2) \, = \, \inf \{a>0 : \sup_A (\nu_1(A)-\nu_2(A^a)) \leq a \, \textrm{ where }
A^a=\{x : d(x,A)\leq a\} \} \, ,$$  $$d_{FM}(\nu_1,\nu_2) \, = \, \sup \left\{\int f (d\nu_1-d\nu_2) \,
\textrm{ for } f\in BLip(E) \textrm{ such that } \parallel f\parallel_{BLip}\leq 1\right\} \, ,$$ where
$BLip$ is the set of bounded and Lipschitz functions and $\parallel f\parallel_{BLip}=\parallel
f\parallel_\infty + \parallel f \parallel_{Lip}$ . For both metrics $M_1(E)$ is Polish. If in
addition $E$ is compact then so does $M_1(E)$.

Furthermore the following inequalities are known to hold $$d_{FM}(\nu_1,\nu_2)\leq \parallel
\nu_1-\nu_2\parallel_{TV} \quad \textrm{ and } \quad d_{P}(\nu_1,\nu_2)\leq \frac 12 \, \parallel
\nu_1-\nu_2\parallel_{TV} \, ,$$ and $$\phi(d_{P}(\nu_1,\nu_2))\leq d_{FM}(\nu_1,\nu_2) \leq 2
d_{P}(\nu_1,\nu_2) \, ,$$ where $\phi(u)=\frac{2u^2}{2+u}$.
\end{proposition}

In the sequel $$C_n \, = \, C^{\varepsilon_n} \, = \, \{\nu : \bar{d}(C,\nu)\leq \varepsilon_n\}$$
where $\bar{d}$ is one of the previous metrics.

\begin{definition}\label{4.4}
Let $(X,d)$ a metric space. If $A\subseteq X$ is totally bounded, we denote by
$N_X(A,d,\varepsilon)$ the minimal number of (open) balls with radius $\varepsilon$ that covers
$A$. The function $N_X$ is often called the metric entropy. In the sequel we simply note
$N(d,\varepsilon)$ the quantity $N_X(X,d,\varepsilon)$, if $X$ is totally bounded.
\end{definition}

Our first result is concerned with compact state spaces.

\begin{theorem}\label{4.5}
Assume that $E$ is compact. Let $C$ be a narrowly closed convex subset of $M_1(E)$ such that
$H(C\mid \alpha)<+\infty$, and $\alpha^*$ be the $I$- projection of $\alpha$ on $C$. Then for any
sequence $\varepsilon_n$ going to 0 and such that $N(d_{FM},\varepsilon_n/4) \, e^{- \, n
\varepsilon_n^2/8} \, \to \, 0$ (resp. $N(d_{P},\varepsilon_n/4) \, e^{- \, n \varepsilon_n^2/2} \,
\to \, 0$) as $n \to \infty$, $\alpha_{C_n,k}^n \, \to \, \alpha^{* \otimes k}$ in total variation
distance.
\end{theorem}

\begin{proof}
Let $B(\alpha^*,\varepsilon)$ the open ball centered at $\alpha^*$ with radius $\varepsilon$. Then
$$\alpha^{* \otimes n}(L_n \in C^{\varepsilon}) \geq \alpha^{* \otimes n}(L_n \in
B(\alpha^*,\varepsilon)) = \, 1 \, - \, \alpha^{* \otimes n}(L_n \in B^c(\alpha^*,\varepsilon))$$
where $B^c$ is as usual the complement subset of $B$. But we can recover
$B^c(\alpha^*,\varepsilon)$ by $$N_{M_1(E)}(B^c(\alpha^*,\varepsilon),\bar{d},\eta) \leq
N(\bar{d},\eta)$$ closed balls with radius $\eta$ so that $$\alpha^{* \otimes n}(L_n \in
B^c(\alpha^*,\varepsilon)) \, \leq \, N(\bar{d},\eta) \, \max_j \alpha^{* \otimes n}(L_n \in B_j)$$
for such balls $B_j$. But a closed ball being closed and convex, Lemma \ref{2.11} shows that
$$\alpha^{* \otimes n}(L_n \in B_j) \, \leq \, e^{- \, n H(B_j\mid \alpha^*)} \, .$$ Since $B_j \subseteq
(B^c(\alpha^*,\varepsilon))^{2\eta} $ we have $H(B_j\mid \alpha^*) \geq
H((B^c(\alpha^*,\varepsilon))^{2\eta}\mid \alpha^*)$ and finally
$$\alpha^{* \otimes n}(L_n \in B^c(\alpha^*,\varepsilon)) \, \leq \, N(\bar{d},\eta) \,
 e^{- \, n H((B^c(\alpha^*,\varepsilon))^{2\eta}\mid \alpha^*)} \, .$$ Choosing $\eta=\varepsilon/4$,
 hence $(B^c(\alpha^*,\varepsilon))^{2\eta}=B^c(\alpha^*,\varepsilon/2)$, we may apply the results
 recalled in Proposition \ref{4.3} to  get that for all $\nu
 \in B^c(\alpha^*,\varepsilon/2)$, $$H(\nu\mid \alpha^*)\geq \frac 12 \parallel \nu -
 \alpha^*\parallel_{TV}^2 \geq \frac 12 \, d_{FM}^2(\nu,\alpha^*) \geq \varepsilon^2/8 \, .$$ We
 can replace 8 by 2 when replacing the Fortet Mourier metric by the Prohorov one.

 Hence we may apply Corollary \ref{4.2} and Theorem \ref{2.7}.
 \end{proof}

The condition on $\varepsilon_n$ in the previous Theorem is interesting if it can be satisfied by
at least one such sequence. The following proposition shows that it is always the case, it also
relies the metric entropy on $M_1(E)$ to the metric entropy on $E$.

\begin{proposition}\eqlabel{\thetheorem}\label{4.6}
Let $(E,d)$ be a compact metric space. Then for all $\varepsilon >0$,
\begin{enumerate}
\item[(\thetheorem.1)] \quad $N(d_P,\varepsilon) \, \leq \,
\left(\frac{2e}{\varepsilon}\right)^{N(d,\varepsilon)}$ , \item[(\thetheorem.2)] \quad
$N(d_{FM},\varepsilon) \, \leq \, \left(\frac{4e}{\varepsilon}\right)^{N(d,\varepsilon/2)}$,
\item[(\thetheorem.3)] there exists at least one sequence $\varepsilon_n$ going to 0 and such that
$$\lim_{n \to \infty} \left(\frac{n \varepsilon_n^2}{8} + (\log \varepsilon_n) \,
N(d,\varepsilon_n/8)\right) = +\infty \, .$$ Such a sequence fulfills the condition in Theorem
\ref{4.5} for both metrics on $M_1(E)$ (but is not sharp).
\end{enumerate}
\end{proposition}

\begin{proof}
The first result is due do Kulkarni-Zeitouni (\cite{KZ95} Lemma 1), the second one follows thanks
to Proposition \ref{4.3}.

Consider $$f : ]0,1] \to \R^+ \, , \, \varepsilon \, \rightarrow \, - \, \frac{8 (\log \varepsilon)
\, N(d,\varepsilon/8)}{\varepsilon^2} \, ,$$ which is clearly decreasing with infinite limit at 0.
Let $u_n$ a $]0,1]$ valued non increasing sequence, $w_n=f(u_n)$ is then non decreasing with
infinite limit. Introduce for $n$ large enough $k_n=\max \{k \in \N^* \, , \, s.t. w_k\leq
\sqrt{n}\}$.
\begin{itemize}
\item If for all $n$ large enough, $k_n\leq n$, we choose $\varepsilon_n=u_{k_n}$ for all $n\in
[k_n,k_{n+p_n}[$ where $p_n=\inf\{p\geq 1 \, , \, k_{n+p}>k_n\}$. On one hand $n\varepsilon_n^2
\geq k_n \, u^2_{k_n}$ goes to infinity. On the other hand, $$n\varepsilon_n^2 + 8(\log
\varepsilon_n) N(d,\varepsilon_n/8)=n\varepsilon_n^2\left(1 - \frac{w_{k_n}}{n}\right) \geq
n\varepsilon_n^2\left(1 - \frac{1}{\sqrt n}\right) \to +\infty \, .$$ \item If not, there exists
some sequence $p_j$ growing to infinity such that $k_{p_j}\geq p_j$, i.e. $w_{p_j}\leq \sqrt{p_j}$.
Define $\phi(n)$ as the unique integer number such that $n\in [p_{\phi(n)},p_{\phi(n)+1}[$, and
choose $\varepsilon_n=u_{p_{\phi(n)}}$. Then $n\varepsilon_n^2 \geq p_{\phi(n)} \,
u^2_{p_{\phi(n)}}$ goes to infinity and $$n\varepsilon_n^2 + 8(\log \varepsilon_n)
N(d,\varepsilon_n/8)=n\varepsilon_n^2\left(1 - \frac{w_{p_{\phi(n)}}}{n}\right) \geq
n\varepsilon_n^2\left(1 - \frac{1}{\sqrt {p_{\phi(n)}}}\right) \to +\infty \, .$$
\end{itemize}
The final statement is a consequence of the previous ones. The proof is thus completed.
\end{proof}

\begin{example}
If $E$ is a $q$ dimensional compact riemanian manifold, it is known that $N(d,\varepsilon)\leq
C_E/\varepsilon^q$ for some constant $C_E$. In this case we may thus choose $\varepsilon_n = 1/n^a$
for all $0<a<\frac{1}{q+2}$. The size of enlargement is thus much greater than for $F$-moment
constraints.
\end{example}
\medskip

When $E$ is no more compact, but still Polish, it can be  approximated by compact subsets with
large probability. Here are the results in this direction

\begin{theorem}\label{4.8}
Let $C$ be a narrowly closed convex subset of $M_1(E)$ such that $H(C\mid \alpha)<+\infty$, and
$\alpha^*$ be the $I$- projection of $\alpha$ on $C$. Assume that there exist a sequence $(K_n)_n$
of compact subsets of $E$ and a sequence $(\eta_n)_n$ of non negative real numbers such that
$$n\eta_n^2 + 8(\log \eta_n) N_E(K_n,d,\eta_n/8) \, \to \, +\infty$$ as $n \to \infty$. Let
$\varepsilon_n \, = \, \eta_n + 2 \alpha^*(K_n^c)$. If one of the following additional assumptions
\begin{itemize}
\item $\lim_{n \to \infty} \, \left(\alpha^*(K_n)\right)^n \, = \, 1$ , \item $\log
\frac{d\alpha^*}{d\alpha}$ is continuous and bounded, and $\lim_{n \to \infty} \, \alpha^*(K_n) \,
= \, 1$ .
\end{itemize}
Then $\alpha_{C_n,k}^n \, \to \, \alpha^{* \otimes k}$ in total variation distance.
\end{theorem}

Here again the conditions are not sharp, but they hold for both the Prohorov and the Fortet Mourier
metrics.

\begin{proof}
The proof lies on the following Lemma
\begin{lemma}\label{4.9}
For all compact subset $K$ and all $\eta>0$, $$\alpha^{* \otimes n}\left(\bar{d}(L_n,C)\leq \eta +
2 \alpha^*(K^c)\right) \, \geq \, (\alpha^*(K))^n \, \left(1 \, - \, (16 e/\eta)^{N_E(K,d,\eta/8)}
\, e^{- \, n \eta^2/8}\right) \, .$$
\end{lemma}
\begin{proof}{\textit {of the Lemma.}}
Introduce $\alpha^*_K \, = \, \frac{\BBone_K}{\alpha^*(K)} \, \alpha^*$ . Then
\begin{eqnarray*}
\bar{d}(\alpha^*_K,\alpha^*) &\leq & \parallel \alpha^*_K - \alpha^*\parallel_{TV} \, = \, \int \,
\left|\frac{\BBone_K}{\alpha^*(K)} - 1\right| \, d\alpha^* \, \leq \, 2 \, \alpha^*(K^c) \, ,
\end{eqnarray*}
so that according to the triangle inequality $\bar{d}(\nu,\alpha^*) \leq \bar{d}(\alpha^*_K,\nu) +
2\alpha^*(K^c)$ for all $\nu$. Hence $B(\alpha_K^*,\eta) \subseteq \{\nu \, , \, \bar{d}(\nu,C)
\leq \eta + 2 \alpha^*(K^c)\}$ and
\begin{eqnarray*}
\alpha^{* \otimes n}\left(\bar{d}(L_n,C)\leq \eta + 2 \alpha^*(K^c)\right) & \geq & \alpha^{*
\otimes n} (L_n \in B(\alpha_K^*,\eta)) \\ & \geq & \alpha^{* \otimes n} (L_n \in
B(\alpha_K^*,\eta) \textrm{ and } x\in K^n) \\& \geq & (\alpha^*(K))^n \, \alpha_K^{* \otimes n}
(L_n \in B(\alpha_K^*,\eta)) \, .
\end{eqnarray*}
As in the proof of Theorem \ref{4.5} and using (\ref{4.6}.1 or 2) we have $$\alpha_K^{* \otimes n}
(L_n \in B(\alpha_K^*,\eta)) \, \geq \, 1 \, - \, N_{M_1(K)}(\bar{d},\eta/4) \, e^{- \, n\eta^2/8}
\, \geq \, 1 \, - \, (16 e/\eta)^{N_K(d,\eta/8)} \, e^{- \, n\eta^2/8} \, .$$
\end{proof}

The first part of the Theorem is then immediate.
\smallskip

The second part is a little bit more tricky. Let $h=\log \frac{d\alpha^*}{d\alpha}$ . For all
$\varepsilon > 0$
\begin{eqnarray*}
\alpha^{\otimes n}(L_n\in C^{\varepsilon})&\geq& \alpha^{\otimes n}(L_n\in B(\alpha^*,\varepsilon))
\,  = \, \int \BBone_{B(\alpha^*,\varepsilon)}(L_n) \, e^{-n\langle L_n,h\rangle} \, d\alpha^{* \otimes n} \\ &
\geq & e^{-nH(C\mid \alpha)} \, \int \BBone_{B(\alpha^*,\varepsilon)}(L_n) \, e^{-n\langle L_n - \alpha^*,h\rangle }
\, d\alpha^{* \otimes n} \\ & \geq & \ e^{-nH(C\mid \alpha)} \, e^{-n \Delta(\varepsilon)} \, \alpha^{*
\otimes n}(L_n\in B(\alpha^*,\varepsilon))
\end{eqnarray*}
where $\Delta(\varepsilon)=\sup_{\nu \in B(\alpha^*,\varepsilon)} \langle\nu - \alpha^*,h\rangle$ . Since $h$
is continuous and bounded, it is immediate that $\Delta(\varepsilon)$ goes to 0 as $\varepsilon$
goes to 0. Hence if $\varepsilon_n$ goes to 0 $$\liminf_{n \to \infty} \log \left(\alpha^{\otimes
n}(L_n \in C_n) \, e^{n H(C\mid \alpha)}\right) \, \geq \, \liminf_{n \to \infty} \log \left(\alpha^{*
\otimes n}(L_n \in B(\alpha^*,\varepsilon_n))\right) \, .$$ Thus if we choose $\varepsilon_n$ as in
the statement of the Theorem, the right hand side of the previous inequality is greater than
$$\liminf \left(\log \alpha^* (K_n) \, + \, \frac 1n \, \log \left(1 - (16
e/\eta_n)^{N_E(K_n,d,\eta_n/8)} \, e^{- \, n \eta_n^2/8}\right)\right) \, = \, 0$$ and we may apply
Theorem \ref{2.7}.
\end{proof}
In the next section we shall study some typical examples.

\section{\bf Examples.}\label{Exemples}

In Section \ref{Moment} we already discussed the examples of $F$-moments. In this section we shall
first look at the finite dimensional situation, then study examples in relation with Stochastic
Mechanics.

\subsection{\bf Finite dimensional convex constraints.}\label{Exemples1}

\begin{proposition}\label{5.1}
If $E=\R^q$, let $C$ be a narrowly closed convex subset of $M_1(E)$ such that
$H(C\mid \alpha)<+\infty$, and $\alpha^*$ be the $I$- projection of $\alpha$ on $C$. Then
$\alpha_{C_n,k}^n \, \to \, \alpha^{* \otimes k}$ in total variation distance with
$\varepsilon_n=2/n^b$ and $0 < b  <  \frac{1-\frac qa}{2+q}$ provided there exists $a>q$ such that
$\int \parallel x\parallel^a \, d\alpha^* \, < +\infty$ (that holds in particular if $\int \,
e^{\lambda
\parallel x\parallel^a} d\alpha <+\infty$ for some $\lambda>0$).

In addition if either $\int e^{\lambda \parallel x\parallel^a} d\alpha^* <+\infty$ for some
$\lambda>0$, or $\log \frac{d\alpha^*}{d\alpha}$ is bounded and continuous, we may choose $b  <
\frac{1}{2+q}$.
\end{proposition}
Of course in general hypotheses on $\alpha^*$ are difficult to check directly. That is why the
$\alpha$ exponential integrability is a pleasant sufficient condition.

\begin{proof}
Let $M=\int \parallel x\parallel^a \, d\alpha^*$. For $K_n=B(0,n^u)$ we have $$(\alpha^*(K_n))^n \,
\geq \, \left(1 \, - \, \frac{M}{n^{au}}\right)^n \, \to \, 1 \, ,$$ provided $au>1$. In addition
$$N_E(K_n,d,\eta/48) \, \leq \, M' n^{uq}/\eta^q$$ so that if $\eta_n = 1/n^b$ with $b>0$
$$n\eta_n^2 + 8(\log \eta_n) N_E(K_n,d,\eta_n/8) \, \geq \, n^{1-2b} \, \left( 1 \, - \, 8b M' (\log n)
n^{uq+b(2+q)-1}\right)$$ goes to $+\infty$ as soon as $b<\frac{1-uq}{2+q}$ i.e. if $b<\frac{1-\frac
qa}{2+q}$ since $au>1$. We may thus apply Theorem \ref{4.8} with
$\varepsilon_n=(1/n^b)+2(M/n^{ua})\leq 2(1/n^b)$ for $n$ large enough.

If the $\alpha^*$ exponential integrability condition is satisfied we may choose $a$ as large as we
want. If $\log \frac{d\alpha^*}{d\alpha}$ is bounded, $\alpha^*(K_n)$ growing to 1, the condition
$ua>1$ is non necessary.
\end{proof}

\subsection{\bf Schr\"{o}dinger bridges.}\label{Exemples2}

In this subsection and the next one $E=C^0([0,1],M)$ where $M$ is either $\R^q$ or a smooth
connected and compact riemannian manifold of dimension $q$. $E$ is equipped with the sup-norm and
for simplicity with the Wiener measure $\mathcal W$ (i.e. the infinitesimal generator is the
Laplace Beltrami operator), with initial measure $\mu_0$.

An old question by Schr\"{o}dinger can be described as following (see \cite{Fo88} for the original
sentence in french). Let $(X_j)_{j=1,\ldots,n}$ be a $n$-sample of $\mathcal W$. Assume that the
empirical measure at time $1$ (i.e. $L_n(1)=\frac 1n \sum_{j=1}^n \delta_{X_j(1)}$) is far from the
expected law $\mu_1$ of the Brownian Motion at time 1. What is the most likely way to observe such
a deviation ? Clearly the answer (when the number of Brownian particles grows to infinity) is
furnished by the Gibbs conditional principle : the most likely way is to imagine that any block of
$k$ particles is made of (almost) independent particles with common law $\mathcal W^*$ which
minimizes $H(\mathcal V\mid \mathcal W)$ among all probability measures on $E$ such that $\mathcal
W\circ X^{-1}(0)=\mu_0$ and $\mathcal W\circ X^{-1}(1)$ belongs to the observed set of measures. If
the observed set is reduced to a single measure (thin) a double limit formulation of this principle
is contained in the first chapter of \cite{Aebi}.

To be precise introduce for $\varepsilon \geq 0$
\begin{equation}\label{5.2}
C^{\varepsilon}(\nu_0,\nu_1) \, = \, \{\mathcal V \in M_1(E) \, s.t. \, \bar{d}(\mathcal
V_0,\nu_0)\leq \varepsilon \, , \, \bar{d}(\mathcal V_1,\nu_1)\leq \varepsilon\} \, ,
\end{equation}
where $\mathcal V_t$ denotes the law $\mathcal V \circ X^{-1}(t)$. When $\varepsilon =0$ we will
not write the superscript $0$. We are in the situation studied in the previous section since
$C(\nu_0,\nu_1)$ is a narrowly closed convex subset of $M_1(E)$. We shall write $\mathcal W^*$ the
$I$- projection of $\mathcal W$ on $C$ (without specifying unless necessary the initial and final
measures) when it exists.

Before to apply the results in Section \ref{Convex} we shall recall some known results about $C$
and $\mathcal W^*$.

Denote by $\mathcal V_{u,v}$ (resp. $\mathcal W_{u,v}$) the conditional law of $\mathcal V$ knowing
that $X(0)=u$ and $X(1)=v$, i.e. the law of the $\mathcal V$ bridge from $u$ to $v$. Also denote by
$\nu_{0,1}$ (resp. $\mu_{0,1}$) the $\mathcal V$ (resp. $\mathcal W$) joint law of $X(0),X(1)$. The
decomposition of entropy formula $$H(\mathcal V\mid \mathcal W) \, = \, H(\nu_{0,1}\mid \mu_{0,1}) + \int
H(\mathcal V_{u,v}\mid \mathcal W_{u,v}) \, d\nu_{0,1}(u,v) \, ,$$ immediately shows that, if it
exists, $$\mathcal W^* \, = \, \int \mathcal W_{u,v} \, d\mu^*_{0,1}(u,v) \, ,$$ where
$\mu^*_{0,1}$ is the $I$- projection of $\mu_{0,1}$ on $$\Pi(\nu_0,\nu_1)=\{\beta \in M_1(M\times
M) \, s.t. \, \beta_0=\nu_0 \, , \, \beta_1=\nu_1 \} \, ,$$ if it exists. In other words the
problem reduces to a finite dimensional one, i.e. on $M \times M$. The following Theorem collects
some results we need

\begin{theorem}\label{5.3}
Assume that $H(\nu_0\mid \mu_0)$ and $H(\nu_1\mid \mu_1)$ are both finite and that $p = \log
\frac{d\mu_{0,1}}{d(\mu_0 \otimes \mu_1)} \, \in \mathbb L^1(\nu_0 \otimes \nu_1)$. Then
$H(\Pi(\nu_0,\nu_1)\mid \mu_{0,1})$ is finite.

In addition $\frac{d\mu_{0,1}^*}{d\mu_{0,1}} (u,v)=f(u)g(v)$ for any pair of functions $(f,g)$
satisfying
\begin{equation}\label{5.4}
\left\{\begin{array}{l}\frac{d\nu_0}{d\mu_0}(u)=f(u)\int p(u,v)g(v)d\mu_1(v)\\
 \frac{d\nu_1}{d\mu_1}(v)=g(v)\int p(u,v)f(u)d\mu_0(u)
 \end{array}\right..
\end{equation}
\end{theorem}

The proof is contained in \cite{CG99} Proposition 6.3 and \cite{Fo88} p.161-164.
\medskip

Finally under the assumptions of Theorem \ref{5.3} $$\frac{d\mathcal W^*}{d\mathcal W} = f(X(0)) \,
g(X(1)) \, .$$ We can now state

\begin{theorem}\label{5.5}
Under the assumptions of Theorem \ref{5.3}, $$\mathcal W_{\varepsilon_n,k}^n := \mathcal L
(X_1,\ldots,X_k/L_n \in C^{\varepsilon_n}(\nu_0,\nu_1)) \, \to \, \mathcal W^{* \otimes k}$$ in total
variation distance for all sequence $\varepsilon_n$ going to 0 such that the following holds : for
all sequence $(Y_j)_j$ (resp. $(Z_j)_j$) of i.i.d. random variables with law $\nu_0$ (resp.
$\nu_1$) , $$\lim_{n \to \infty} \mathbb P(\bar d (L_n^Y,\nu_0)\leq \varepsilon_n) = 1 \textrm{ and
} \lim_{n \to \infty} \mathbb P(\bar d(L_n^Z,\nu_1)\leq \varepsilon_n) = 1 \, .$$ In particular the
above convergence holds for instance in the following two cases
\begin{itemize}
\item $M$ is compact and $n\varepsilon_n^2 + 8(\log \varepsilon_n) N_M(d,\varepsilon_n/8) \, \to \,
+\infty$, \item $M=\R^q$ , there exists $a>q$ such that for $i=0,1$ , $\int \parallel x\parallel^a
d\nu_i < +\infty$ , $\varepsilon_n=2/n^b$ and $ b< \frac{1-\frac qa}{2+q}$.
\end{itemize}
\end{theorem}

For the proof just apply Corollary \ref{4.2}, and for the examples Proposition \ref{4.6} and
Proposition \ref{5.1}.

\subsection{\bf Nelson processes.}\label{Exemples3}

A natural generalization of the framework of Subsection \ref{Exemples2} is to impose the full flow
of marginal laws instead of only the initial and final ones. Building diffusion processes with a
given flow of marginal laws is the first step in Nelson's approach of the Schr\"{o}dinger equation.
The problem was first tackled by Carlen \cite{Ca84}. Relationship with minimization of entropy was
first observed by H. F\"{o}llmer (\cite{Fo88}) and explored in details in a series of papers by C.
L\'eonard and the first named author (\cite{CL94,CL95,CL96}). This approach and the results below
can be viewed as some ``statistical mechanics'' approach of quantum mechanics. We shall not discuss
further the meaning of the previous sentence here. We prefer insist on the enormous difference
between a pair and the flow of all marginal laws.

Hence here $$C(\nu_t) \, = \, \{\mathcal V \in M_1(E) \, s.t. \forall t \in [0,1] \, , \, \mathcal
V_t=\nu_t\} \, ,$$ and for $\varepsilon >0$ $$C^{\varepsilon}(\nu_t) \, = \, \{\mathcal V \in
M_1(E) \, s.t. \bar d (\mathcal V,C(\nu_t)) \leq \varepsilon\} \, .$$ For simplicity we shall only
consider the case $M=\R^q$ (though a similar discussion is possible for a general connected and
compact riemannian manifold). Not to lose sight of our main goal we first state the convergence
result we have in mind, and will discuss the hypotheses later on.

\begin{theorem}\label{5.6}
Assume that $C(\nu_t)$ is non empty and that $\mathcal W$ has an $I$- projection $\mathcal W^*$ on
$C(\nu_t)$, such that $\log \frac{d\mathcal W^*}{d\mathcal W}$ is bounded and continuous. Assume in
addition that the initial law $\mu_0$ has a polynomial concentration rate i.e. $\mu_0(B(0,R))\leq
C/R^m$ for some $m>0$ and all $R>0$. Then if $\varepsilon_n=1/(\log n)^r$ for some $r<1/2q$ ,
$$\mathcal W_{\varepsilon_n,k}^n := \mathcal L (X_1,\ldots,X_k/L_n \in C^{\varepsilon_n}(\nu_t)) \, \to
\, \mathcal W^{* \otimes k}$$ in total variation distance.
\end{theorem}

\begin{proof}
According to Theorem \ref{4.8} it is enough to find a sequence $K_n$ of compact subspaces of $E$
and a sequence $\eta_n$ of positive numbers going to 0 such that $$\lim_{n \to \infty} \mathcal
W^*(K_n)=1 \, \textrm{ and } \, \lim_{n \to \infty} \left(n \eta_n^2 + 8 (\log \eta_n) \,
N_E(K_n,\parallel \, \parallel_{\infty},\eta_n/8)\right) = +\infty \, .$$ Since $\frac{d\mathcal
W^*}{d\mathcal W}$ is bounded by some $e^D$ , we may replace the first condition by $\lim_{n \to
\infty} \mathcal W(K_n)=1$ and choose $\varepsilon_n\geq \eta_n + 2 e^D \mathcal W(K_n^c)$. The
most natural way to choose such compact sets is to use Kolmogorov regularity criterion. Since the
support of $\mathcal W$ is included into the set of H\"{o}lder paths of order $\beta<1/2$ introduce
$$K(R,M,\beta) \, = \, \left\{w\in E \, s.t. |w(o)|\leq R \, \textrm{ and } \, \sup_{s\neq t \in
[0,1]}\frac{\|w(s)-w(t)\|}{|s-t|^\beta}\leq M\right\} \, ,$$ for $R$, $M$ positive and $\beta<1/2$.
 Kolmogorov's criterion tells us that $$\mathcal W (K^c(R,M,\beta)) \, \leq \,
 \mu_0(B(0,R))+C(p,\beta) M^{-p}$$ for all $p>1$. In addition, thanks to Theorem 2.7.1 p.155 in
 \cite{VW95} $$N_E(K(R,M,\beta),\parallel \, \parallel_{\infty},\eta/8) \, \leq \, c_1(\beta,q) \,
 (8R/\eta)^q \, e^{c_2(\beta,q) (M/\eta)^{q/\beta}} \, .$$ Choosing $K_n=k(R_n,M_n,\eta_n)$ with
 $$R_n=(a \log n)^{\beta/qm} \quad M_n=(b \log n)^{\beta/q} \quad \eta_n=(c \log n)^{-\beta/q}$$ we
 see that $n \eta_n^2 + 8 (\log \eta_n) \,
N_E(K_n,\parallel \, \parallel_{\infty},\eta_n/8)$ is less than $$n (\log n)^{-2\beta/q} \,
\left(A_1+A_2 \log(c \log n) \, (\log n)^{q+\frac{2\beta}{q}} \, n^{c_2(\beta,q)bc-1}\right)$$ for
some $A_1$ and $A_2$ independent of $n$. Choosing $b$ in such a way that $c_2(\beta,q)bc-1<0$  we
obtain a leading term going to $+ \infty$ as $n$ goes to $\infty$.

Putting all this together, we get $$\eta_n + 2 e^D \mathcal W(K_n^c)\leq (c \log n)^{-\beta/q} + 2
C e^D (a \log n)^{-\beta/q} + 2 e^D C(p,\beta) (b \log n)^{-\beta p/q}$$ which is less than $(\log
n)^{-\beta'/q}$ for all $\beta'<1/2$ and $n$ large enough.
\end{proof}

\begin{remark}\label{5.7}
The assumption $\log \frac{d\mathcal W^*}{d\mathcal W}$  bounded and continuous is essential.
Indeed without it Theorem \ref{4.8} requires $(\mathcal W^*(K_n))^n$ goes to 1, i.e. $\mathcal
W^*(K_n^c)=o(1/n)$. Assuming that $\frac{d\mathcal W^*}{d\mathcal W}$ belongs to $\mathbb
L_{r}(\mathcal W)$, Kolmogorov criterion yields  $M_n$  of order $n^a$. It is then easy to see that
this is no more compatible with any choice of $\eta_n$ such that $\lim_{n \to \infty} \left(n
\eta_n^2 + 8 (\log \eta_n) \, N_E(K_n,\parallel \, \parallel_{\infty},\eta_n/8)\right) = +\infty$.
\end{remark}
\medskip

To conclude this subsection let us say a few words about our assumptions.

First of all $C(\nu_t)$ is non empty as soon as $\nu_t$ satisfies a Fokker-Planck equation with a
drift $B(t,X(t))$ of finite energy (i.e. $\int_0^1 \int B^2(t,x) d\nu_t dt <+\infty$) see
\cite{Ca84,CL94,CL95}. In addition Girsanov theory is still available (see \cite{CL94,CL95} for the
details) so that $$\frac{d\mathcal W^*}{d\mathcal W}=\frac{d\nu_0}{d\mu_0} \, \exp \left(\int_0^T
B(t,w(t)) dw(t) - 1/2 \int_0^T |B(t,w(t))|^2 dt\right)$$ where $T=\inf \{s \leq 1 \, s.t. \int_0^s
|B(t,w(t))|^2 dt = + \infty\}$. In general this density (even when $T=1$) is not continuous.
\smallskip

Nevertheless some interesting cases enter the framework of Theorem \ref{5.6}.

Let $U$ be a $C^2_b$ potential. Then the law $\mathcal V_0$ of the unique strong solution of
$$dX_t=dW_t - \nabla U(X_t) dt \quad , \quad \mathcal L(X_0)=\nu_0$$ satisfies $$\frac{d\mathcal V_0}{d\mathcal W}
=\frac{d\nu_0}{d\mu_0} \, \exp \left(U(w(0))-U(w(1)) - 1/2 \int_0^1 (|\nabla U|^2-\Delta U)(t,w(t))
dt\right) \, .$$ Hence $\log \frac{d\mathcal V_0}{d\mathcal W}$ is bounded and continuous as soon
as $\log \frac{d\nu_0}{d\mu_0}$ is. In addition $\mathcal V_0$ is the $I$- projection of $\mathcal
W$ on $C(\nu_t)$ where $\nu_t=\mathcal L(X_t)$ (see \cite{CL94}). The conclusion of Theorem
\ref{5.6} is thus available for $\mathcal V_0$. If we replace $\R^q$ by a compact manifold we may
include the stationary (actually reversible) case i.e. $\nu_0=e^{-2U} dx/Z_U$.

\section{\bf A super-thin case: volatility calibration.}\label{Calib}

In subsections \ref{Exemples2} and \ref{Exemples3} we have studied the laws of some diffusion
processes from the point of view of $I$- projections, hence we only allowed a change of drift. We
shall now study the opposite situation: the drift being fixed, how to choose the diffusion
coefficient. We thus immediately lose any kind of absolute continuity, introducing a new difficulty
that is super-thin subsets. Let us describe precisely the problem.
\medskip

Consider a family (indexed by continuous time-space functions $\sigma$) of S.D.E.
\begin{equation}\label{6.1}
\forall t \in[0,1],\quad dX(t)=\s(t,X(t)) \, dw(t)+b_0(t,X(t)) \, dt \quad ; \quad X(0)=0 \, ,
\end{equation}
where $w$ is a standard Brownian motion. We assume that $b_0$ is continuous and bounded and $$0<  \s_{min} \leq \s
\leq  \s_{max} <+\infty$$ for some real numbers $\s_{min}$ and $ \s_{max}$. Under this assumption, it is well known that (\ref{6.1}) admits weak solutions and that there is uniqueness in law. We will denote in the sequel $\Q_{\s,\,b_0}$ the probability measure on $\Omega=C([0,1],\R)$ thus defined by (\ref{6.1}).

In \cite{AFHS} the authors addressed the problem of calibrating $\s$ (volatility in mathematical
finance) when $b_0$ is known (a consequence of the ``absence of arbitrage'') and $X$ satisfies a set
of generalized moment constraints
\begin{equation}\label{6.2}
\E \, [f_j(t_j,X(t_j))]=c_j \quad , \quad j\in\Lambda \, , \,  \Lambda \hbox{ finite.}
\end{equation}
Their strategy is based on the following Bayesian principle : take a prior $\s_0$,
the corresponding prior law of $X$ is $\Q_{\s_0,\,b_0}$. Then the ``most probable'' $\P$ satisfying \eqref{6.2}, will be
the one which minimizes the relative entropy $H(\P\mid \Q_{\s_0,\,b_0})$. Of course this principle is meaningless
here. Indeed, the finiteness of $H(\P\mid \Q_{\s_0,\,b_0})$ implies that $\P$ has the same diffusion coefficient
as $\Q_{\s_0,\,b_0}$, hence there is no such $\P$ satisfying \eqref{6.2} unless $\Q_{\s_0,\,b_0}$ does. To bypass this
difficulty, the authors propose to approximate $\Q_{\s_0,\,b_0}$ by some well chosen $\Q_{\s_0,\,b_0}^{\ep}$ (actually
various time discretization), in such a way that $\ep \, H(\P^{\ep}\mid \Q_{\s_0,\,b_0}^{\ep})$ goes to some limit
$K(\P\mid \Q_{\s_0,\,b_0})$, and then use $K$ as the cost function to be minimized.
\smallskip

We shall interpret this strategy in the following way.

For simplicity assume that the set of constraints is reduced to a single one i.e. introduce the set
$$C_F \, = \, \left\{\P \, , \, \E_{\P} \, [F(X(1))]=1 \right\}$$ where $\P$ describes the set of Probability
measures on $\Omega=C([0,1],\R)$. We will choose as before some $\varepsilon$ enlargement of $C_F$,
i.e. define $$C_F^\varepsilon  \, = \, \left\{ \P \, , \, \left|\int \,F(X(1)) \, d\P \, - \,
1\right|<\varepsilon\right\} \, .$$ Again for simplicity, we shall assume that $b(t,x)=b_0$
for some $b_0>0$ (extensions to more general cases can be easily done).
We also define
$$\Sigma_0=\left\{\sigma: [0,1]\times \R \rightarrow ]\sigma_{min},\sigma_{max}[, \text{ continuous}\right\}$$
and for $\varepsilon <b_0$,
$$\mathcal{B}_\varepsilon=\left\{b: [0,1]\times \R\rightarrow ]b_0-\varepsilon,b_0+\varepsilon[, \text{ continuous}\right\}.$$
Let us precise that the space of space-time continuous functions $C([0,1]\times\R,\R)$ will always be furnished with the topology
 of uniform convergence on every compact subset of $[0,1]\times\R$.

Now we introduce a standard approximation of $\Q_{\s,\,b_0}$, namely the trinomial tree.

Choose some $\alpha > \s_{max}$ and $0 < s < b_0$. For $(y,z)\in \R^2$ we define
$$\left\{
\begin{array}{l}
     m^{n}(y,z)=\frac{y^2}{2\alpha^2}+\frac{z}{2\alpha \sqrt{n}} \\
     d^{n}(y,z)=\frac{y^2}{2\alpha^2}-\frac{z}{2\alpha \sqrt{n}} \\
     r^{n}(y,z)=1-\frac{y^2}{\alpha^2}\\
\end{array}
\right..$$ For $n$ large enough ($>n_0$), it is easily seen that for all $(y,z)\in
[\s_{min},\s_{max}]\times [b_0-s,b_0+s]$ the vector $(m^n,d^n,r^n)$ has all its entries strictly
positive (their sum being 1), so that we may define the following transition kernel defined on $\R$
for all $(\sigma,b) \in \overline{\Sigma_0} \times \overline{\mathcal{B}_s}$, $n \geq n_0$ and
$(t,x)\in [0,1]\times \R$,
$$\Pi^n_{\sigma,\,b}(t,x,\,.\,)=m^n(\sigma,b)(t,x).\delta_{x+\frac{\alpha}{\sqrt{n}}}+
r^n(\sigma,b)(t,x).\delta_{x}+d^n(\sigma,b)(t,x).\delta_{x-\frac{\alpha}{\sqrt{n}}}.$$

We thus define the probability measure  $\Q^n_{\sigma,\,b}$
\begin{equation}\label{deftrinom}
\left\{
\begin{array}{ll}
    (1)& \Q^n_{\sigma,\,b}(X_0=0)=1,  \\
    (2)& \Q^n_{\sigma,\,b}\left(X_t=X_{\frac{k}{n}}+(nt-k)\left[X_{\frac{k+1}{n}}-X_{\frac{k}{n}}\right],\frac{k}{n}\leq t \leq \frac{k+1}{n}\right)=1,  \\
    (3)& \Q^n_{\sigma,\,b}\left(X_{\frac{k+1}{n}}\in \,.\,\left|X_{\frac{k}{n}},\ldots ,X_0\right.\right)= \Pi_{\sigma,\,b}^n\left(\frac{k}{n},X_{\frac{k}{n}},\,.\,\right)\\
\end{array}
\right.
\end{equation}\\

In the sequel, we will denote by $\E^n_{\sigma,\,b}[\,.\,]$ the expectation with respect to the trinomial tree $\Q^n_{\s,\,b}$. The support of $\Q^n_{\s,\,b}$ is $\Omega_n\subset \Omega$
defined by
$$\Omega_n=\left\{\omega \in \Omega :
\left[
\begin{array}{ll}
    - \omega(0)=0 &  \\
    - \omega\left(\frac{i+1}{n}\right)-\omega\left(\frac{i}{n}\right)\in
\left\{-\frac{\alpha}{\sqrt{n}},0,\frac{\alpha}{\sqrt{n}}\right\}, & \hbox{pour $i=0,\ldots ,n-1$} \\
    - \omega \text{ affine on} \left[\frac{i}{n},\frac{i+1}{n}\right], & \hbox{pour $i=0,\ldots ,n-1$} \\
\end{array}
\right. \right\}$$ The set $\Omega_n$ is finite with cardinality $3^n$.


\medskip

Finally denoting by $L_m=\frac 1m \, \sum_{i=1}^m \, \delta_{\omega_i}$ the empirical measure on
$\Omega$, we shall study $\R^{n,m}_\varepsilon $ defined by $$\R^{n,m}_\varepsilon (B)\, = \,
(\Q^n_{\s_0,b_0})^{\otimes m}(\omega_1\in B / L_m \in \widetilde{\T}^n_\varepsilon \cap C_F^\varepsilon) \, ,$$
where $\widetilde{\T}^n_\ep$ will be defined later. Let us just say for the moment that $\widetilde{\T}^n_\ep$
is an open set of $M_1(\Omega_n)$ which contains all the trinomial trees $\Q^n_{\s,\,b}$ with $\s$
in a totally bounded subset $\Sigma_1$ of $\Sigma_0$ and $b\in \mathcal{B}_\ep$.
Roughly speaking, for each level of approximation ($n$) we consider a $m$ sample of the trinomial tree and
look at the conditional law of the first coordinate, knowing that the empirical measure is not too
far from being a trinomial tree satisfying the moment constraint.

Our aim is to show that one can find sequences $\varepsilon_n$ going to 0 and $m_n$ going to
infinity, such that $\R^{n,m_n}_{\varepsilon_n}$ goes towards some $\Q_{\s^*,b_0}$, the one
proposed in \cite{AFHS} we will now describe.
\smallskip

First, for fixed $n$ and $\varepsilon$, since all measures are defined on a finite set, it is not
difficult to see that the set $\mathcal{M}^n_\ep$ of minimizers of $H(\,.\,\mid Q^n_{\s_0,\,b_0})$ on $\widetilde{\T}^n_\varepsilon \cap C_F^\varepsilon$ is nonempty. It can then be shown that the elements of $\mathcal{M}^n_\ep$ are still a trinomial trees. Now an easy
computation shows that $\displaystyle{\sigma \mapsto \frac{1}{n} \,
H(\Q^n_{\sigma,\,b}\mid \Q^n_{\sigma_0,\,b_0})}$ is converging (in a sense close to the
$\Gamma$-convergence sense) to $$\sigma\mapsto I(\sigma\mid \sigma_0)=\E_\sigma\left[\int_0^1
q(\sigma^2(X_t,t),\sigma_0^2(t,X_t))\,dt\right],$$ with
$$q(x,y)=\log\left(\frac{x}{y}\right)\frac{x}{\alpha^2}+
\log\left(\frac{\alpha^2-x}{\alpha^2-y}\right)\left[1-\frac{x}{\alpha^2}\right].$$ One thus expects
that the limit $\Q_{\s^*,b_0}$ is the one obtained by minimizing $I$ on $\overline{ \Sigma_0}$ under the moment constraint.
\medskip

The remainder of this section will be devoted to give rigorous statements and proofs. Note that the
result gives a rigorous statistical flavor to the method proposed by Avellaneda et altri.

\subsection{Presentation of the results}
We recall that the space $C([0,1]\times\R,\R)$ is equipped with the topology of uniform
convergence on every compact subsets of $[0,1]\times\R$. Before presenting our results, let us
state the basic convergence property of trinomial trees :
\begin{proposition}\label{cv-Trinom}
If $s\geq\ep_n\geq0$ goes to zero and $\s_n\in \overline{\Sigma_0}$
goes to $\s \in \overline{\Sigma_0}$ then, for all $b_n \in
\overline{\mathcal{B}_{\ep_n}}$, the sequence
     $\Q^n_{\s_n,\,b_n}$ goes to $\Q_{\sigma,\,b_0}$.
\end{proposition}
From now on, we will make the following assumptions :
\begin{itemize}
\item The minimum value of the function $I(\,.\,\mid \s_0)$ on the set
$\left\{\sigma \in \overline{ \Sigma_0}:\int
F(X_1)\,dQ_{\s,\,b}=1\right\}$ is attained at a unique point $\s^*$.
\item The minimizer $\s^*$ belongs to $\Sigma_0$.
\end{itemize}
Now let us introduce some notations. For all $\s\in
\overline{\Sigma_0}$, let $\Delta_{n,\,\s}$ be the continuity modulus
of $\s$ on the compact set $[0,1]\times
\left[-\a\sqrt{n},\a\sqrt{n}\right]$, ie.
\begin{equation*}
\Delta_{n,\,\s}(\varepsilon)=\sup\bigg\{|\s(t,x)-\s(s,y)| : s,t \in
[0,1], x,y \in \left[-\a\sqrt{n},\a\sqrt{n}\right], |t-s|+|x-y|\leq
\varepsilon\bigg\}.
\end{equation*}
Let $\Sigma_1$ be defined by
$$\Sigma_1=\{\s \in \Sigma_0 :\forall n \in \N^*,\quad \Delta_{n,\,\s}< 2\Delta_{n,\,\s^*}\}.$$
According to Ascoli Theorem, $\Sigma_1$ is easily seen to be totally bounded.\\
Now let us consider the set $\widetilde{\T}^n_\ep$ of all probability
measures $\Q$ on $\Omega$ satisfying
\begin{equation}\label{tildeT}
\left\{
\begin{array}{ll}
    (1)& \Q(X_0=0)=1,  \\
    (2)& \Q\left(X_t=X_{\frac{k}{n}}+(nt-k)\left[X_{\frac{k+1}{n}}-X_{\frac{k}{n}}\right],\frac{k}{n}\leq t \leq \frac{k+1}{n}\right)=1,  \\
    (3)& \exists (\s,b) \in \Sigma_1\times \mathcal{B_\varepsilon} \text{ such that }
\Q\left(X_\frac{p+1}{n}\in
\,.\,\left|X_\frac{p}{n}\right.\right)=\Pi^n_{\s,\,b}\left(\frac{p}{n},X_\frac{p}{n},\,.\,\right)
\end{array}
\right.
\end{equation}
In the sequel we will set $A^n_{\ep}:=\widetilde{\T}^n_\ep\cap C_F^\ep$. Defining (when possible),
for all positive integer $m$, $$ \R^n_{\ep,\,m}=\E_{(\Q_{\sigma_0,\,b_0}^n)^{\otimes
m}}\left[L_m\left|L_m\in A^n_{\ep}\right.\right],
$$
our main result is the following :

\begin{theorem}\label{cali-mr}
If $\ep^0_n=\min\left(\left|\E^n_{\s^*,\,b_0}\left[
F\left(X_1\right)\right]-1\right|+1/n, s\right)$, then there exists
a sequence $m_n$ of positive integers going to $+\infty$, such that
$\R^n_{\ep^0_n,\,m_n}$ converges to $\Q_{\s^*,\,b_0}$.
\end{theorem}

In order to prove this theorem, the first step is to study the
convergence of $\R^n_{\ep^0_n,m}$ when $n$ is fixed and $m$ goes to
$+\infty$. This is done in the two following propositions :

\begin{proposition}\label{Gibbs-nfixed}
Recall that $d_{FM}$ denotes the Fortet-Mourier distance, and for
all $\ep>0$ let $\mathcal{M}^n_\ep$ be the set of minimizers of
$H(.\mid Q^n_{\s_0,\,b_0})$ on $\overline{A^n_\ep}$. Then,
$$d_{FM}(\R^n_{\ep^0_n,m},\overline{co}\, \mathcal{M}^n_{\ep^0_n})\xrightarrow[m\rightarrow+ \infty]{}0,$$
where $\overline{co}\, \mathcal{M}^n_{\ep^0_n}$ denotes the closed convex hull of $
\mathcal{M}^n_{\ep^0_n}$.
\end{proposition}
\proof The set $A^n_{\ep^0_n}$ is non empty (it contains $\Q^n_{\s^*,\,b_0}$) and, according to the
 proposition below, it is open and satisfies $H(A^n_{\ep^0_n}\mid
Q^n_{\s_0,\,b_0})=H(\overline{A^n_{\ep^0_n}}\mid Q^n_{\s_0,\,b_0})$ . The result follows
immediately from the classical Gibbs conditioning principle.
\endproof
\begin{proposition}\label{prop-cond-event}~
\begin{enumerate}
\item The set $A^n_\ep$ is an open subset of $M_1(\Omega_n)$, and satisfies
$H(A^n_\ep\mid Q^n_{\s_0,\,b_0})=H(\overline{A^n_\ep}\mid Q^n_{\s_0,\,b_0})$.
\item  Every element
of $\mathcal{M}^n_\ep$ is of the form $\Q^n_{\s,\,b}$ for some
$(\s,b) \in \overline{Z_1}\times \overline{\mathcal{B}_\ep}$.
\end{enumerate}
\end{proposition}
According to Proposition \ref{Gibbs-nfixed}, we know that for large $m$, $\R^n_{\ep^0_n,m}$ is
close to $\overline{co}\, \mathcal{M}^n_{\ep^0_n}$. The next step consists in proving that this set
is close to $\{\Q_{\s^*,\,b_0}\}$. This will follow from the particular type of convergence of the
normalized entropy functions :
\begin{proposition}\label{gammaconv}~
\begin{enumerate}
    \item If $0<\ep_n$ goes to $0$,
     then for every sequence $b_n \in \mathcal{B}_{\ep_n}$, and for every $\s \in \overline{\Sigma_0}$,
     the following holds :
$$\frac{H(\Q^n_{\s,\,b_n}\mid \Q^n_{\s_0,\,b_0})}{n}\xrightarrow[n \rightarrow +\infty]{} I(\s\mid \s_0).$$
\item Furthermore, if $\s_n \in \overline{\Sigma_0}$ converges to $\s \in \overline{\Sigma_0}$, then
$$\liminf_{n \rightarrow +\infty}\frac{H(\Q^n_{\s_n,\,b_n}\mid \Q^n_{\s_0,\,b_0})}{n}\geq I(\s\mid \s_0).$$
\end{enumerate}
\end{proposition}
\begin{remark}Recall that a sequence $f_n$ of real valued functions
defined on some metric space $\Gamma$-converges to some function $f$, if
\begin{itemize}
\item for all $x$, $\lim_{n\rightarrow +\infty}f_n(x)=f(x)$, \item for all sequence $x_n$
converging to some $x$, $\liminf_{n\rightarrow +\infty} f_n(x_n)\geq f(x)$.
\end{itemize}
The preceding proposition can thus be restated by saying that for
every $b_n\in \mathcal{B}_{\ep_n}$ with $\ep_n$ going to $0$, the
sequence of functions $\s \mapsto
\frac{H(\Q^n_{\s,\,b_n}\mid \Q^n_{\s_0,\,b_0})}{n}$ $\Gamma$-converges
to $\s \mapsto I(\s\mid \s_0)$.\\
It is well known that this kind of convergence is well adapted for
deriving the convergence of minimizers. The next proposition
illustrates this fact :
\end{remark}
\begin{proposition}\label{conv-min}
Suppose that for every $n$, $\Q^n_{\s_n,\,b_n}$ is an element of
$\mathcal{M}^n_{\ep^0_n}$, then
\begin{equation}\label{eq:conv-min}
\Q^n_{\s_n,\,b_n}\xrightarrow[n\rightarrow +
\infty]{}\Q_{\s^*,\,b_0}.
\end{equation}
\end{proposition}
\proof For all $n$, $\Q^n_{\s^*,\,b_0}$ belongs to $A^n_{\ep^0_n}$.
Thus, using the minimization property of $\Q^n_{\s_n,\,b_n}$, one
has $\frac{1}{n}H(\Q^n_{\s_n,\,b_n}\mid \Q^n_{\s_0,\,b_0})\leq
\frac{1}{n}H(\Q^n_{\s^*,\,b_0}\mid \Q^n_{\s_0,\,b_0})$. According to
point (1) of Proposition \ref{gammaconv}, this implies that
\begin{equation}\label{eq:conv-min1}
\limsup_{n\rightarrow +
\infty}\frac{1}{n}H(\Q^n_{\s_n,\,b_n}\mid \Q^n_{\s_0,\,b_0})\leq
I(\s^*\mid \s_0).
\end{equation}
According to the point (2) of Proposition \ref{prop-cond-event}, $\s_n \in \overline{\Sigma_1}$.
This set being compact, one can find some converging subsequence $\s_{n_p}$. Let $\tilde{\s}$ be
its limit. The point (2) of Proposition \ref{gammaconv}, yields :
\begin{equation}\label{eq:conv-min2}
\liminf_{p\rightarrow + \infty}
\frac{1}{n_p}H(\Q^{n_p}_{\s_{n_p},\,b_{n_p}}\mid \Q^{n_p}_{\s_0,\,b_0})\geq
I(\tilde{\s}\mid \s_0).
\end{equation}
From (\ref{eq:conv-min1}) and (\ref{eq:conv-min2}), one deduces that
$$I(\tilde{\s}\mid \s_0)\leq I(\s^*\mid \s_0).$$

As $\s^*$ is the unique minimizer of $I(.\mid \s_0)$ under the moment constraint, one has
$\tilde{\s}=\s^*$. The point $\s^*$ is thus the unique accumulation point of the compact sequence
$\s_n$. It follows that $\s_n$ converges to $\s^*$. Now, (\ref{eq:conv-min}) follows immediately
from Proposition \ref{cv-Trinom}.
\endproof
We are now ready to prove Theorem \ref{cali-mr}.\\
\emph{Proof of Theorem \ref{cali-mr}.} First, we have the following
immediate inequality
$$d_{FM}\left(\R^n_{\ep^0_n,\,m},\Q_{\s^*,\,b_0}\right)
\leq d_{FM}\left(\R^n_{\ep_n,\,m},\overline{co}\,\mathcal{M}^n_{\ep^0_n}\right)+ \sup_{\Q \in\,
\overline{co}\,\mathcal{M}^n_{\ep^0_n}}d_{FM}\left(\Q,\Q_{\s^*,\,b_0}\right).$$ Thus, according to
Proposition \ref{Gibbs-nfixed}, it suffices to prove that $$\sup_{\Q \in
\,\overline{co}\,\mathcal{M}^n_{\ep^0_n}}d_{FM}\left(\Q,\Q_{\s^*,\,b_0}\right)\xrightarrow[n\rightarrow+\infty]{}0.$$
The application $\Q \mapsto d_{FM}\left(\Q,\Q_{\s^*,\,b_0}\right)$ being convex and continuous, we
get $$\sup_{\Q \in
\,\overline{co}\,\mathcal{M}^n_{\ep^0_n}}d_{FM}\left(\Q,\Q_{\s^*,\,b_0}\right)=\sup_{\Q \in
\,\mathcal{M}^n_{\ep^0_n}}d_{FM}\left(\Q,\Q_{\s^*,\,b_0}\right).$$ But $\mathcal{M}^n_{\ep_n^0}$ is
compact. Thus, there exists $\Q^n_{\s_n,\,b_n}\in \mathcal{M}^n_{\ep^0_n}$, such that
$$\sup_{\Q \in \,\mathcal{M}^n_{\ep^0_n}}d_{FM}\left(\Q,\Q_{\s^*,\,b_0}\right)=
d_{FM}\left(\Q^n_{\s_n,\,b_n},\Q_{\s^*,\,b_0}\right).$$
Applying Proposition \ref{conv-min}, we get
$$\Q^n_{\s_n,\,b_n}\xrightarrow[n \rightarrow + \infty]{} \Q_{\sigma^*,\,b_0},$$
which achieves the proof. \hfill $\square$\\
Before giving the proofs of Proposition
\ref{prop-cond-event} and \ref{gammaconv}, let us do some comments
on our result.
\begin{remark}~
\begin{itemize}
\item The reason why we work with $\widetilde{\T}^n_\ep$ instead of the more natural set
$\T^n_\ep=\{\Q^n_{\s,\,b} : \s \in \Sigma_1, b\in \mathcal{B}_\ep\}$ is that $\T^n_\ep$ is of
\emph{empty interior}. The set $\T^n_\ep$ was thus a bad candidate for defining a conditioning
event in Gibbs Principle. In fact, from the relative entropy point of view, working with
$\widetilde{\T}^n_\ep$ does not change anything : point (2) of Proposition \ref{prop-cond-event}
shows that the entropy minimizers on $A^n_\ep$ \emph{are trinomial trees}. \item We introduced the
set $\Sigma_1$ because some compactness is needed in Proposition \ref{prop-cond-event}. Note that
if we replace $\Sigma_1$ by $\Sigma_0$ in the definition of $\widetilde{\T}^n_\ep$, this set
becomes \emph{convex} (see \cite{G}). In this framework, there is a unique entropy-minimizer
$\Q^n_{\sigma_n^*,\,b_n^*}$. But we are not able to prove directly that the sequence $\s_n^*$ is
compact. If this was true, Theorem \ref{cali-mr} would hold with $\Sigma_0$ replacing $\Sigma_1$.
\item The assumption that $I(.\mid \s_0)$ admits a unique minimizer under the moment constraint is
needed in the proof of Theorem \ref{cali-mr}. Namely, we used in the proof the fact that the
function $\Q\mapsto d_{FM}(\Q,\Q_{\s^*,\,b_0})$ is convex. If we were dealing with a set
$\mathcal{M}$ of minimizers containing more than one element, this function would be replaced by
the function $\Q\mapsto d_{FM}(\Q,\mathcal{M})$ which is no longer convex.
\end{itemize}
\end{remark}

\bigskip

\subsection{Proofs}~\\
\emph{Proof of (1) of Proposition \ref{prop-cond-event}.}
The set $C_F^ep$ being clearly open, it suffices to show that $\widetilde{\T}^n_\ep$ is an open subset of $M_1(\Omega_n)$.
First, it is easily seen that there is a constant $c>0$ depending only on $\s_{min}$, $\s_{max}$, $b_0$, $s$ and $\a$ such that
$$\Q\left(X_{\frac{k}{n}}=\frac{j\alpha}{\sqrt{n}}\right)> c,$$
for all $\Q \in \widetilde{\T}^n_\ep$ and all $|j|\leq k \leq n$.
For all $|j|\leq k \leq n$ and $\Q \in M_1(\Omega_n)$, let us define
\begin{equation}\label{Fkj}
F_{k,\,j}(\Q)=\a \sqrt{n}\frac{\Q\left(X_{\frac{k+1}{n}}=\frac{(j+1)\a}{\sqrt{n}},X_{\frac{k}{n}}=\frac{j\a}{\sqrt{n}}\right)-
\Q\left(X_{\frac{k+1}{n}}=\frac{(j-1)\a}{\sqrt{n}},X_{\frac{k}{n}}=\frac{j\a}{\sqrt{n}}\right)}{\Q\left(X_{\frac{k}{n}}=\frac{j\a}{\sqrt{n}}\right)}
\end{equation}
and
\begin{equation}\label{Gkj}
G_{k,\,j}(\Q)=\a^2 \frac{\Q\left(X_{\frac{k+1}{n}}=\frac{(j+1)\a}{\sqrt{n}},X_{\frac{k}{n}}=\frac{j\a}{\sqrt{n}}\right)+
\Q\left(X_{\frac{k+1}{n}}=\frac{(j-1)\a}{\sqrt{n}},X_{\frac{k}{n}}=\frac{j\a}{\sqrt{n}}\right)}{\Q\left(X_{\frac{k}{n}}=\frac{j\a}{\sqrt{n}}\right)}
\end{equation}
These applications are continuous on the open set $$\left\{\Q\in M_1\left(\Omega_n\right) : \forall |j|\leq k \leq n,\quad \Q\left(X_{\frac{k}{n}}=\frac{j\a}{\sqrt{n}}\right)> c\right\}$$ and the following holds
$$\Q\in \widetilde{\T}^n_\ep\Leftrightarrow \left\{\begin{array}{l}
\forall |j|\leq k \leq n\\
\forall |q|\leq p \leq n \\
\end{array},\quad \begin{array}{l}
\Q\left(X_{\frac{k}{n}}=\frac{j\a}{\sqrt{n}}\right)> c,\\
F_{k,\,j}(Q)\in]b_0-\ep,b_0+\ep[,\\
G_{k,\,j}(Q)\in]\s^2_{min}, \s^2_{max}[,\\
\left|\sqrt{G_{k,\,j}(\Q)}-\sqrt{G_{p,\,q}(\Q)}\right|<2\Delta_{n,\,\sigma^*}\left(\left|\frac{k}{n}-\frac{p}{n}\right|+\left|\frac{\a j}{\sqrt{n}}-\frac{\a q}{\sqrt{n}}\right|\right)
\end{array}\right.$$
One easily concludes from this that $\widetilde{\T}^n_\ep$ is an open subset of $M_1(\Omega_n)$.\\
Now let us show that $H(A^n_\ep\mid Q^n_{\s_0,\,b_0})=H(\overline{A^n_\ep}\mid Q^n_{\s_0,\,b_0})$. As $\Q^n_{\s_0,\,b_0}$ gives a positive mass to every trajectory of $\Omega_n$, the convex function $M_1(\Omega_n)\ni \Q \mapsto H(\Q\mid \Q^n_{\s_0,\,b_0})$ is everywhere finite thus continuous. As a consequence, $H(O\mid \Q^n_{\s_0,\,b_0})=H(\overline{O}\mid \Q^n_{\s_0,\,b_0})$ holds true for all open set $O$ of $M_1(\Omega_n)$. This is in particular true for $A^n_\ep$.\hfill $\square$\\

In order to prove the point (2) of Proposition \ref{prop-cond-event}, we need the following lemma.\\
\begin{lemma}\label{lemmaH}
For all $\s\in \overline{\Sigma_0}$, $b\in \mathcal{B}_\ep$, $\ep\leq s$, let us define :
$$q^n_{\s,\,b\,;\,\s_0,\,b_0}(t,x,y)=\cfrac[l]{d\Pi_{\s,\,b}^n(t,x,\,.\,)}{d\Pi_{\s_0,\,b_0}^n(t,x,\,.\,)}(y)$$
and
$$h^n_{\s,\,b\,;\,\s_0,\,b_0}(t,x)=H(\Pi_{\s,\,b}^n(t,x,\,.\,)\mid \Pi_{\s_0,\,b_0}^n(t,x,\,.\,)).$$
Then it holds :
\begin{align}
\cfrac[l]{d\Q_{\s,\,b}^n}{d\Q_{\s_0,\,b_0}^n}&=\prod_{i=0}^{n-1}q^n_{\s,\,b\,;\,\s_0,\,b_0}
\left(\frac{i}{n},X_\frac{i}{n},X_\frac{i+1}{n}\right)\label{lemmaH1}\\
H(\Q_{\s,\,b}^n\mid \Q_{\s_0,\,b_0}^n)&=\sum_{i=0}^{n-1}\E_{\s,\,b}^n\left[h^n_{\s,\,b\,;\,\s_0,\,b_0}\left(\frac{i}{n},X_\frac{i}{n}\right)\right]\label{lemmaH2}
\end{align}
Let $\Q$ be a probability measure satisfying
\begin{equation}\label{tildeTbis}
\left\{
\begin{array}{ll}
    (1)& \Q(X_0=0)=1,  \\
    (2)& \Q\left(X_t=X_{\frac{k}{n}}+(nt-k)\left[X_{\frac{k+1}{n}}-X_{\frac{k}{n}}\right],\frac{k}{n}\leq t \leq \frac{k+1}{n}\right)=1,  \\
    (3)& \Q\left(X_\frac{p+1}{n}\in
\,.\,\left|X_\frac{p}{n}\right.\right)=\Pi^n_{\s,\,b}\left(\frac{p}{n},X_\frac{p}{n},\,.\,\right)
\end{array}
\right.,
\end{equation}
for some $\s\in \overline{\Sigma_0}$ and $b\in \overline{\mathcal{B}_\ep}$.
Then
\begin{equation}\label{lemmaH3}
\forall i=0,\ldots,n-1,\quad \mathcal{L}_\Q\left(X_\frac{i}{n}\right)=\mathcal{L}_{\Q_{\sigma,\,b}^n}\left(X_\frac{i}{n}\right).
\end{equation}
Furthermore,
\begin{equation}\label{lemmaH4}
H(\Q\mid \Q_{\sigma_0,\,b_0}^n)=H(\Q\mid \Q_{\sigma,\,b}^n)+H(\Q_{\sigma,\,b}^n\mid \Q_{\sigma_0,\,b_0}^n).
\end{equation}
\end{lemma}
\proof The proofs of (\ref{lemmaH1}), (\ref{lemmaH2}), (\ref{lemmaH3}) rely on very easy
computations and are left to the reader. Let us prove (\ref{lemmaH4}). It is clear that,
\begin{equation}\label{lemmaH5}
H(\Q\mid \Q_{\s_0,\,b_0}^n)  =H(\Q\mid \Q_{\s,\,b}^n)+\int\log\left(\cfrac[l]{d\Q_{\s,\,b}^n}{d\Q_{\s_0,\,b_0}^n}\right)d\Q.
\end{equation}
Next, we have
\begin{align*}
&\int\log\left(\cfrac[l]{d\Q_{\sigma,\,b}^n}{d\Q_{\sigma_0,\,b_0}^n}\right)d\Q  \stackrel{(i)}{=}  \int \sum_{i=0}^{n-1}\log\left[q^n_{\sigma,\,b\,;\,\sigma_0,\,b_0}\left(\frac{i}{n},X_\frac{i}{n},X_\frac{i+1}{n}\right)\right]d\Q\\
& \stackrel{(ii)}{=} \E_\Q\left[\sum_{i=0}^{n-1}\int\log\left[q^n_{\s,\,b\,;\,\s_0,\,b_0}\left(\frac{i}{n},X_\frac{i}{n},y\right)\right]\Pi_{\s,\,b}^n\left(\frac{i}{n},X_\frac{i}{n},dy\right)\right]\\
& =  \sum_{i=0}^{n-1}\E_\Q\left[h^n_{\s,\,b\,;\,\s_0,\,b_0}\left(\frac{i}{n},X_\frac{i}{n}\right)\right]
 \stackrel{(iii)}{=}  \sum_{i=0}^{n-1}\E_{\s,\,b}^n\left[h^n_{\s,\,b\,;\,\s_0,\,b_0}\left(\frac{i}{n},X_\frac{i}{n}\right)\right]\\
& \stackrel{(iv)}{=}  H(\Q_{\s,\,b}^n\mid \Q_{\s_0,\,b_0}^n),
\end{align*}
where (i) follows from (\ref{lemmaH1}), (ii) is obtained by conditioning by $X_i$, (iii) is a consequence of (\ref{lemmaH3}) and (iv) of (\ref{lemmaH2}).
 Plugging this in (\ref{lemmaH5}), we obtain (\ref{lemmaH4}).
\endproof

\emph{Proof of (2) of Proposition \ref{prop-cond-event}.} Let $\Q$ be in $\mathcal{M}^n_\ep$. As
$\Q$ belongs to $\overline{A^n_\ep}$, there exist $\s \in \overline{\Sigma_1}$ and $b\in
\overline{\mathcal{B}_\ep}$ such that (\ref{tildeTbis}) is fulfilled. According to (\ref{lemmaH4}),
one has
$$H(\Q\mid \Q_{\s_0,\,b_0}^n)=
H(\Q\mid \Q_{\s,\,b}^n)
+H(\Q_{\s,\,b}^n\mid \Q_{\s_0,\,b_0}^n).$$
If $\Q_{\s,\,b}^n$ belongs to $\overline{A^n_\ep}$, then we deduce from the preceding equation
 that $H(\overline{A^n_\ep}\mid \Q^n_{\s_0,\,b_0})\geq H(\Q\mid \Q_{\s,\,b}^n)+ H(\overline{A^n_\ep}\mid \Q^n_{\s_0,\,b_0})$,
 and consequently $H(\Q\mid \Q_{\s,\,b}^n)=0$, which implies that $\Q=\Q^n_{\s,\,b}$. Thus, the only thing to do is to prove
  that $\Q_{\s,\,b}^n\in\overline{A^n_\ep}$.

Let $\left(\Q_p\right)_p$ be a sequence of $A^n_\ep$ going to $\Q$. For each $p$, there is a pair
$(\s_p,b_p)\in \Sigma_1\times\mathcal{B}_\ep$ such that (\ref{tildeTbis}) is fulfilled. For all
$|j|\leq k\leq n$, one has
$$b_p\left(\frac{k}{n},\frac{\a j}{\sqrt{n}}\right)=F_{k,\,j}(\Q_p)$$
and
$$\s_p^2\left(\frac{k}{n},\frac{\a j}{\sqrt{n}}\right)=G_{k,\,j}(\Q_p),$$
where $F_{k,\,j}$ and $G_{k,\,j}$ are defined by (\ref{Fkj}) and (\ref{Gkj}). These functions being continuous, we have
$$b_p\left(\frac{k}{n},\frac{\a j}{\sqrt{n}}\right)\xrightarrow[p \rightarrow +\infty]{} b\left(\frac{k}{n},\frac{\a j}{\sqrt{n}}\right)$$
and
$$\s_p^2\left(\frac{k}{n},\frac{\a j}{\sqrt{n}}\right)\xrightarrow[p \rightarrow +\infty]{}\left(\s\right)^2\left(\frac{k}{n},\frac{\a j}{\sqrt{n}}\right),$$
for all $|j|\leq k\leq n$. It follows easily that
$$\Q^n_{\s_p,\,b_p}\xrightarrow[p\rightarrow+\infty]{}\Q_{\s,\,b}^n.$$
But according to (\ref{lemmaH3}),
$$\Q_p \in A^n_\ep \Rightarrow \Q^n_{\s_p,\,b_p}\in A^n_\ep.$$
Consequently, $\Q_{\s,\,b}^n$ is in the closure of $A^n_\ep$.\hfill $\square$\\

\emph{Proof of Proposition \ref{gammaconv}.} Recall that for all
$\sigma \in \overline{\Sigma_0}$, $I(\s\mid \s_0)$ is defined by
\begin{align*}
I(\s\mid \s_0)&=\E_{\s,\,b}\left[\int_0^1q(\s^2(t,X_t),\s^2_0(t,X_t))\,dt\right],\\
\intertext{with}
q(x,y)&=\log\left(\frac{x}{y}\right)\frac{x}{\a^2}+\log\left(\frac{\a^2-x}{\a^2-y}\right)\left[1-\frac{x}{\a^2}\right].
\end{align*}

(1) Let us show that there exists some $K>0$, depending only on $\a$, $\s_{min}$, $\s_{max} $,
$b_0$ and $s$, such that
\begin{equation}\label{DL}
\left|h^n_{\s,\,b\,;\,\s_0,\,b_0}-q(\s^2,\s^2_0)\right|\left(\frac{k}{n},x\right)
\leq \frac{K}{n}
\end{equation}
for all $(k,x) \in \{0,\ldots ,n-1\}\times\frac{\a}{\sqrt{n}}\mathbb{Z}$ and $(\s,b)\in \overline{\Sigma_0}\times\overline{\mathcal{B}_s}$.\\
For all $(\s,b)\in \overline{\Sigma_0}\times\overline{\mathcal{B}_s}$ :
\begin{align*}
\log\left[\cfrac[l]{m^n(\s,b)}{m^n(\s_0,b_0)}\right]m^n(\s,b)&=
\left[\log\left(\!\frac{\s^2}{\s_0^2}\right)\!+\log\left(\!1+\frac{b\a}{\sqrt{n}\s^2}\right)\!-
\log\left(\!1+\frac{b_0\a}{\sqrt{n}\s_0^2}\right)\right]
 \times\left[\frac{\s^2}{2\a^2}+\frac{b}{2\a\sqrt{n}}\right]\\
\log\left[\cfrac[l]{d^n(\s,b)}{d^n(\s_0,b_0)}\right]d^n(\s,b)&=
\left[\log\left(\!\frac{\s^2}{\s_0^2}\right)\!+\log\left(\!1-\frac{b\a}{\sqrt{n}\s^2}\right)\!-
\log\left(\!1-\frac{b_0\a}{\sqrt{n}\s_0^2}\right)\right]
\times\left[\frac{\s^2}{2\a^2}-\frac{b}{2\a \sqrt{n}}\right]\\
\log\left[\cfrac[l]{r^n(\s,b)}{r^n(\s_0,b_0)}\right]r^n(\s,b) &=
\log\left(\frac{\a^2-\s^2}{\a^2-\s_0^2}\right)\left[1-\frac{\s^2}{\a^2}\right]
\end{align*}

Using Taylor's formula, it is easily seen that for $\varepsilon \in\{-1,1\}$,
$$\sup_{\substack{x \in[\s_{min}^2,\s_{max}^2]\\ y\in[b_0-s,b_0+s]}}
\left|\log\left(1+\frac{\ep y\a}{\sqrt{n}x}\right)-\frac{\ep
y\a}{\sqrt{n}x}+\frac{1}{2}\left(\frac{\ep
y\a}{\sqrt{n}x}\right)^2\right|\leq \frac{K}{n\sqrt{n}},$$ with $K$
depending only on $\a$, $\s_{max} $,
$\s_{min}$, $b_0$ et $s$.\\
After some easy computations, one derives (\ref{DL}) from these inequalities.\\

In the sequel we will use the following notations
$$\Phi^n=\frac{1}{n}\sum_{i=0}^{n-1}q\left(\s^2\left(\frac{i}{n},X_{\frac{i}{n}}\right),\s_0^2\left(\frac{i}{n},X_{\frac{i}{n}}\right)\right)$$
and
$$\Phi=\int_0^1 q(\s^2(t,X_t),\s_0^2(t,X_t))\,dt.$$

The function $q$ is bounded and continuous  on $[\sigma_{min}^2,\sigma_{max}^2]^2$. $\Phi^n$ is
thus a sequence of uniformly bounded continuous functions on $\Omega$, which converges pointwise to
the bounded continuous function $\Phi$. Let us show that $\Phi^n$ converges uniformly to $\Phi$ on
every compact subset of $\Omega$. The function $q$ is Lipschitz on
$[\sigma_{min}^2,\sigma_{max}^2]^2$ ; let $M>0$ be such that
$$|q(x,y)-q(x',y')|\leq M(|x-x'|+|y-y'|).$$
Let $\Delta$ be the continuity modulus of $\s^2$, ie.
$$\Delta(u)=\sup_{|t-s|+|y-x|\leq u}|\s^2(s,x)-\s^2(t,y)|,$$
and $\Delta_0$ the continuity modulus of $\s_0^2$.\\
With these notations, we have
\begin{align*}
\left|\Phi^n-\Phi\right|& =  \left|\frac{1}{n}\sum_{i=0}^{n-1}q\left(\sigma^2\left(\frac{i}{n},X_{\frac{i}{n}}\right),\sigma_0^2\left(\frac{i}{n},X_{\frac{i}{n}}\right)\right)-
\int_0^1 q(\sigma^2(t,X_t),\sigma_0^2(t,X_t))\,dt
\right|\\
& \leq \sum_{i=0}^{n-1} \int_{\frac{i}{n}}^{\frac{i+1}{n}}\left|q\left(\sigma^2\left(\frac{i}{n},X_{\frac{i}{n}}\right),\sigma_0^2\left(\frac{i}{n},X_{\frac{i}{n}}\right)\right)-q(\sigma^2(t,X_t),\sigma_0^2(t,X_t))\right|dt\\
& \leq  M\sum_{i=0}^{n-1} \int_{\frac{i}{n}}^{\frac{i+1}{n}}\left|\sigma^2\left(\frac{i}{n},X_{\frac{i}{n}}\right)-\sigma^2(t,X_t)\right|+
\left|\sigma_0^2\left(\frac{i}{n},X_{\frac{i}{n}}\right)-\sigma_0^2(t,X_t)\right|dt\\
& \leq  M\left[\sup_{|s-t|\leq\frac{1}{n}}\left|\sigma^2(s,X_s)-\sigma^2(t,X_t)\right|+
\sup_{|s-t|\leq\frac{1}{n}}\left|\sigma_0^2(s,X_s)-\sigma_0^2(t,X_t)\right|\right]\\
& \leq  M\left[\sup_{|s-t|\leq\frac{1}{n}}\Delta\left(|s-t|+|X_s-X_t|\right)+\sup_{|s-t|\leq\frac{1}{n}}\Delta_0\left(|s-t|+|X_s-X_t|\right)\right]\\
& \leq  M\left[\Delta\left(\frac{1}{n}+\sup_{|s-t|\leq\frac{1}{n}}|X_s-X_t|\right)+\Delta_0\left(\frac{1}{n}+\sup_{|s-t|\leq\frac{1}{n}}|X_s-X_t|\right)\right]
\end{align*}
Let $\mathcal{K}$ be a compact subset of $\Omega$. According to
Ascoli Theorem, we have
$$\sup_{\omega \in\, \mathcal{K}}\sup_{|t-s|\leq\frac{1}{n}}|X_s-X_t|\xrightarrow[n\rightarrow+\infty]{}0.$$
Thus
$$\sup_{\omega \in\, \mathcal{K}}\left|\Phi^n(\omega)-\Phi(\omega)\right|\xrightarrow[n\rightarrow+\infty]{}0.$$

According to (\ref{DL}) :
$$\left| \frac{1}{n}H(\Q_{\s,\,b_n}^n\mid \Q_{\s_0,\,b_0}^n)-\E_{\s,\,b_n}^n\left[\Phi^n\right]\right|\leq \frac{K}{n}$$
where $K$ depends only on $\a$, $\s_{\text{max}} $,
$\s_{min}$, $b_0$ and $s$. Using the uniform convergence of $(\Phi^n)_n$ on every compact and the tightness of the sequence $\Q_{\s,\,b_n}^n$, it is now easy to see that

$$\lim_{n \rightarrow \infty}\frac{1}{n}H(\Q_{\s,\,b_n}^n\mid \Q_{\s_0,\,b_0}^n)=I(\s\mid \s_0).$$

(2)
\begin{align*}
\frac{1}{n}H(\Q_{\s_n,\,b_n}^n\mid \Q_{\s_0,\,b_0}^n) & =  \frac{1}{n}\int \log\left(\cfrac[l]{d\Q_{\s_n,\,b_n}^n}{d\Q_{\s_0,\,b_0}^n}\right)d\Q_{\s_n,\,b_n}^n\\
& =  \frac{1}{n}\int\!\! \log\!\left(\cfrac[l]{d\Q_{\s_n,\,b_n}^n}{d\Q_{\s,\,b_n}^n}\right)\!\!d\Q_{\s_n,\,b_n}^n\!+\!\frac{1}{n}\int\!\! \log\!\left(\cfrac[l]{d\Q_{\s,\,b_n}^n}{d\Q_{\s_0,\,b_0}^n}\right)\!\!d\Q_{\s_n,\,b_n}^n\\
& =  \frac{1}{n}H(Q_{\s_n,\,b_n}^n\mid  Q_{\s,\,b_n}^n)+\frac{1}{n}\int \log\left(\cfrac[l]{d\Q_{\s,\,b_n}^n}{d\Q_{\s_0,\,b_0}^n}\right)d\Q_{\s_n,\,b_n}^n\\
& \geq  \frac{1}{n}\int \log\left(\cfrac[l]{d\Q_{\s,\,b_n}^n}{d\Q_{\s_0,\,b_0}^n}\right)d\Q_{\s_n,\,b_n}^n\\
\end{align*}

According to (\ref{lemmaH1}) of  Lemma \ref{lemmaH}
$$\frac{1}{n}\int \log\left(\cfrac[l]{d\Q_{\sigma,\,b_n}^n}{d\Q_{\sigma_0,\,b_0}^n}\right)d\Q_{\sigma_n,\,b_n}^n=\E_{\sigma_n,\,b_n}^n\left[\frac{1}{n}\sum_{i=0}^{n-1}k^n\left(\frac{i}{n},X_{\frac{i}{n}}\right)\right],$$
where
$$k^n=\log\left(\cfrac[l]{m^n(\s,b_n)}{m^n(\s_0,b_0)}\right)m^n(\s_n,b_n)+\log\left(\cfrac[l]{r^n(\s,b_n)}{r^n(\s_0,b_0)}\right)r^n(\s_n,b_n)
+\log\left(\cfrac[l]{d^n(\s,b_n)}{d^n(\s_0,b_0)}\right)d^n(\s_n,b_n)$$

It is easily seen that there is a constant $K$ depending only on $\a$, $\s_{min}$,
$\s_{max}$, $b_0$ and $s$ such that
$$\forall R>0,\quad\sup_{|x|\leq R,\, t\in[0,1]}|k^n-h^n_{\s,\,b_n\,;\,\s_0,\,b_0}|(t,x)\leq K\sup_{|x|\leq R,\, t\in[0,1]}|\s_n-\s|(t,x).$$
The sequence $\Q^n_{\s_n,\,b_n}$ converging to
$\Q_{\s,\,b}$, it is a tight sequence. As a consequence, for all $\b>0$, there is $R>0$ such that
$$\Q^n_{\s_n,\,b_n}\left(\sup_{t\in[0,1]}\left|X_t\right|\leq R\right)\geq 1-\beta.$$
One can find $M>0$ depending on $\alpha$, $\s_{min}$,
$\s_{max}$, $b_0$ and $s$, such that $\left|k^n\right|\leq M$ and
$\left|h^n_{\s,\,b_n\,;\,\s_0,\,b_0}\right|\leq M$. Thus,
\begin{align*}
&\left|\E_{\s_n,\,b_n}^n\left[\frac{1}{n}\sum_{i=0}^{n-1}k^n\left(\frac{i}{n},X_{\frac{i}{n}}\right)\right]-
\E_{\s_n,\,b_n}^n\left[\frac{1}{n}\sum_{i=0}^{n-1}h^n_{\s,\,b_n\,;\,\s_0,\,b_0}\left(\frac{i}{n},X_{\frac{i}{n}}\right)\right]\right|\\
&\leq \E_{\s_n,\,b_n}^n\left[\frac{1}{n}\sum_{i=0}^{n-1}\left|h^n_{\s,\,b_n\,;\,\s_0,\,b_0}-k^n\right|\BBone_{[0,R]}(\sup_{t\in[0,1]}|X_t|)\right]
+2M(1-\beta)\\
&\leq K\sup_{|x|\leq R,\,
t\in[0,1]}|\sigma_n-\sigma|(t,x)+2M(1-\beta).
\end{align*}
One easily concludes that
$$\E_{\sigma_n,\,b_n}^n\left[\frac{1}{n}\sum_{i=0}^{n-1}k^n\left(\frac{i}{n},X_{\frac{i}{n}}\right)\right]-
\E_{\sigma_n,\,b_n}^n\left[\frac{1}{n}\sum_{i=0}^{n-1}h^n_{\sigma,\,b_n\,;\,\sigma_0,\,b_0}\left(\frac{i}{n},X_{\frac{i}{n}}\right)\right]\xrightarrow[n\rightarrow
+\infty]{} 0.$$ A similar reasoning as in the proof of point (1) shows that
$$\E_{\sigma_n,\,b_n}^n\left[\frac{1}{n}\sum_{i=0}^{n-1}h^n_{\sigma,\,b_n\,;\,\sigma_0,\,b_0}\left(\frac{i}{n},X_{\frac{i}{n}}\right)\right]\xrightarrow[n \rightarrow +\infty]{} I(\sigma\mid \sigma_0),$$
which achieves the proof. \hfill $\square$

\bibliographystyle{plain}

\end{document}